\documentclass{gtart_h}  

\def\ifplaintex{\expandafter\ifx\csname documentclass\endcsname\relax}

\ifplaintex 
\else
\fi

\expandafter\ifx\csname beginpicture\endcsname\relax
\expandafter\ifx\csname documentclass\endcsname\relax
\input pictex \else
\input prepictex \input pictex \input postpictex \fi\fi

\def\gt{{\mathsurround=0pt\it $\cal G\mskip-2mu$eometry \&\ 
$\cal T\!\!$opology}}        

\def\gtp{{\mathsurround=0pt\it $\cal G\mskip-2mu$eometry \&\ 
$\cal T\!\!$opology $\cal P\!$ublications}}  

\def\lognumber#1{\def\thelognumber{#1}}
\def\volumenumber#1{\def\thevolumenumber{#1}}
\def\papernumber#1{\def\thepapernumber{#1}}
\def\volumeyear#1{\def\thevolumeyear{#1}}

\def\pagenumbers#1#2{\def\startpage{#1}\def\finishpage{#2}}
\def\published#1{\def\publishdate{#1}}
\def\proposed#1{\def\theproposer{#1}}
\def\seconded#1{\def\theseconders{#1}}
\def\received#1{\def\receiveddate{#1}}
\def\revised#1{\def\reviseddate{#1}}
\def\accepted#1{\def\accepteddate{#1}}

\long\def\asciiabstract#1{\long\def\theasciiabstract{#1}}

\let\\\par\let\thelognumber\relax
\let\thevolumenumber\relax\let\thepapernumber\relax
\let\thevolumeyear\relax\let\thesamplenumber\relax\let\startpage\relax
\let\finishpage\relax\let\publishdate\relax\let\receiveddate\relax
\let\reviseddate\relax\let\accepteddate\relax\let\theasciititle\relax
\let\theasciiauthors\relax
\let\theasciiabstract\relax
\let\theasciiemail\relax\let\theshortauthors\relax\let\theshorttitle\relax

\long\def\maketitlep{   

\count0=\startpage

\gt\hfill      
\beginpicture
\setcoordinatesystem units <0.33truein, 0.33truein> point at 2.2 0.9
\setplotsymbol ({$\cal G$})
\plotsymbolspacing=9truept
\circulararc 315 degrees from 0 1 center at 0 0
\setplotsymbol ({$\cal T$})
\circulararc 315 degrees from 1 -1 center at 1 0
\endpicture
\break
{\small\ifx\thesamplenumber\relax 
Volume \else Sample
\fi\thevolumenumber\ (\thevolumeyear)
\startpage--\finishpage\nl
Published: \publishdate}
\vglue 0.5truein plus 0.4fil minus 0.1truein

{\parskip=0pt\leftskip 0pt plus 1fil\def\\{\par\smallskip}{\ifplaintex\large
\else\Large\fi\bf\thetitle}\par\medskip}   

\vglue 0pt plus 0.1fil 

{\parskip=0pt\leftskip 0pt plus 1fil\def\\{\par}{\sc\theauthors}
\par\medskip}

\vglue 0pt plus 0.1fil 

{\small\parskip=0pt\let\newline\\
{\leftskip 0pt plus 1fil\def\\{\par}{\sl\theaddress}\par}
\expandafter\ifx\theemail\relax    
\relax\else\vglue 5pt plus 0.02fil minus 2pt\def\\{\stdspace{\rm 
and}\stdspace} 
\cl{Email:\stdspace\tt\theemail}\fi
\ifx\theurl\relax                  
\relax\else\vglue 5pt plus 0.02fil minus 2pt\def\\{\stdspace{\rm 
and}\stdspace}
\cl{URL:\stdspace\tt\theurl}\fi\par}

\vglue 7pt plus 0.3fil minus 3pt

{\bf Abstract}
\vglue 5pt plus 0.1fil minus 2pt

\theabstract

\vglue 7pt plus 0.3fil minus 3pt

{\bf AMS Classification numbers}\quad Primary:\quad \theprimaryclass

Secondary:\quad \thesecondaryclass

\vglue 5pt plus 0.3fil minus 2pt

{\bf Keywords:}\quad \thekeywords

\vglue 10pt plus 0.5fil minus 5pt

{\small  Proposed: \theproposer\hfill Received: \receiveddate\nl
Seconded: \theseconders\hfill 
\ifx\reviseddate\relax                         
Accepted: \accepteddate                        
\else
Revised: \reviseddate                          
\fi}
\eject
}       

\let\maketitlepage\maketitlep
\let\maketitle\maketitlepage

\font\phead=cmsl9 scaled 950
\font=cmsl9 scaled 1050
\font\pnum=cmbx10 scaled 913
\font\lnum=cmbx10 
\font\pfoot=cmsl9 scaled 950
\font=cmsl9 scaled 1050
\ifplaintex
\headline{\vbox to 0pt{\vskip -4.5mm\line{\small\phead\ifnum
\count0=\startpage ISSN 1364-0380 (on line)
1465-3060 (printed) \hfill {\pnum\folio}\else\ifodd\count0\def\\{ }%
\ifx\theshorttitle\relax\thetitle\else\theshorttitle\fi\hfill{\pnum\folio}
\else\def\\{ and }{\pnum\folio}\hfill\ifx\theshortauthors\relax\theauthors
\else\theshortauthors\fi\fi\fi}\vss}}
\footline{\vbox to 0pt{\vglue 0mm\line{\small\pfoot\ifnum\count0=\startpage
\copyright\ \gtp\hfill\else
\gt, Volume \thevolumenumber\ (\thevolumeyear)\hfill\fi}\vss
}}
\else
\makeatletter
\def\@oddhead{{\small\lhead\ifnum\count0=\startpage ISSN 1364-0380 (on line)
1465-3060 (printed) \hfill {\lnum\number\count0}\else\ifodd\count0
\def\\{ }\ifx\theshorttitle\relax \thetitle \else\theshorttitle\fi\hfill
{\lnum\number\count0}\else\def\\{ and }{\lnum\number\count0}
\hfill\ifx\theshortauthors\relax 
\theauthors\else\theshortauthors\fi\fi\fi}}\def\@evenhead{\@oddhead}
\def\@oddfoot{\small\lfoot\ifnum\count0=\startpage\copyright\ \gtp\hfill\else
\gt, Volume \thevolumenumber\ (\thevolumeyear)\hfill\fi}
\def\@evenfoot{\@oddfoot}
\makeatother
\fi

\newwrite\gtoutfile
\long\gdef\makeheadfile{  
{\def\\{, }\def\s{ }
\immediate\openout\gtoutfile head.xxx
\immediate\write\gtoutfile{Proxy-for: \ifx\theasciiauthors\relax
\theauthors\else\theasciiauthors\fi\s<\ifx\theasciiemail\relax\theemail\else\theasciiemail\fi>}
\immediate\write\gtoutfile{\noexpand\\}
\immediate\write\gtoutfile{Authors: \ifx\theasciiauthors\relax
\theauthors\else\theasciiauthors\fi}
{\def\\{ }\immediate\write\gtoutfile{Title: \ifx\theasciititle\relax
\thetitle\else\theasciititle\fi}}
\immediate\write\gtoutfile{Subj-class: GT or SG or MG etc}
\immediate\write\gtoutfile{MSC-class: \theprimaryclass\ifx\thesecondaryclass\relax\else, \thesecondaryclass\fi}
\immediate\write\gtoutfile{Journal-ref: Geom. Topol. \thevolumenumber
(\thevolumeyear) \startpage-\finishpage}
\immediate\write\gtoutfile{Comments: Published by Geometry and Topology at}
\immediate\write\gtoutfile{\s\s http://www.maths.warwick.ac.uk/gt/GTVol\thevolumenumber/paper\thepapernumber.abs.html}
\immediate\write\gtoutfile{\noexpand\\}
\immediate\write\gtoutfile{}
\ifx\theasciiabstract\relax
\immediate\write\gtoutfile{\theabstract}\else
\immediate\write\gtoutfile{\theasciiabstract}\fi
\immediate\write\gtoutfile{}
\immediate\write\gtoutfile{\noexpand\\}
\immediate\write\gtoutfile{}
\immediate\closeout\gtoutfile}}  

\def\maketitlepage{\maketitlep\makeheadfile}
\let\maketitle\maketitlepage

\lognumber{334}

\volumenumber{8}\papernumber{22}\volumeyear{2004}
\pagenumbers{831}{876}
\received{30 June 2003}
\published{2 June 2004}
\revised{1 June 2004}
\accepted{23 May 2004}
\proposed{Cameron Gordon}
\seconded{Joan Birman, Wolfgang Metzler}

\usepackage{amsmath, amssymb, epsf}

\theoremstyle{definition}
\newtheorem{dfn}{Definition}[section]
\newtheorem{defn}[dfn]{Definition}
\newtheorem{rem}[dfn]{Remark}

\newtheorem{ex}[dfn]{Example}

\theoremstyle{plain}
\newtheorem{thm}[dfn]{Theorem}
\newtheorem{lem}[dfn]{Lemma}

\newtheorem{prop}[dfn]{Proposition}

\def\surface{{\it surface}}

\def\>{\rangle}
\def\<{\langle}

\def\3{\ss}
\def\8{\infty}

\begin{document}

\title{Heegaard splittings of graph manifolds}
\author{Jennifer Schultens}
\address{Department of Mathematics\\1 Shields 
Avenue\\University of California\\Davis, CA 95616, USA}
\email{jcs@math.ucdavis.edu}

\begin{abstract}
Let $M$ be a totally orientable graph manifold with characteristic
submanifold $\cal T$ and let $M = V \cup_S W$ be a Heegaard splitting.
We prove that $S$ is standard.  In particular, $S$ is the amalgamation
of strongly irreducible Heegaard splittings.  The splitting surfaces
$S_i$ of these strongly irreducible Heegaard splittings have the
property that for each vertex manifold $N$ of $M$, $S_i \cap N$ is
either horizontal, pseudohorizontal, vertical or pseudovertical.
\end{abstract}

\asciiabstract{%
Let M be a totally orientable graph manifold with characteristic
submanifold T and let M = V cup_S W be a Heegaard splitting.
We prove that S is standard.  In particular, S is the amalgamation
of strongly irreducible Heegaard splittings.  The splitting surfaces
S_i of these strongly irreducible Heegaard splittings have the
property that for each vertex manifold N of M, S_i cap N is
either horizontal, pseudohorizontal, vertical or pseudovertical.}

\maketitlepage

\section{Introduction}

The subject of this investigation is the structure of Heegaard
splittings of graph manifolds.  This investigation continues the work
begun in \cite{Sch}, \cite{Sch1}, \cite{MS} and \cite{Sch3}.  Since
the publication of those papers, new techniques have been added to the
repertoire of those interested in describing the structure of Heegaard
splittings.  These include the idea of untelescoping a weakly
reducible Heegaard splitting into a generalized strongly irreducible
Heegaard splitting due to M Scharlemann and A Thompson.  They also
include the Rubinstein--Scharlemann graphic, as employed by D Cooper
and M Scharlemann in \cite{CS}.  These insights have not left the
investigation here unaffected.  We hope that their role here is a
tribute to proper {\it affinage}\footnote{This is a French noun
describing the maturing process of a cheese.}.  (The structural
theorem given here has been promised for rather a long time.)  A
similar theorem was announced by J\,H Rubinstein.

The main theorems are the following, for defintions see Sections 2, 3,
4 and 5:

\begin{thm} \label{thm:premain}
Let $M$ be a totally orientable generalized graph manifold.  If $M = V
\cup_S W$ is a strongly irreducible Heegaard splitting, then $S$ is
standard.

More specifically, $S$ can be isotoped so that for each vertex
manifold $M_v$ of $M$, $S \cap M_v$ is either horizontal,
pseudohorizontal, vertical or pseudovertical and such that for each
edge manifold $M_e$, $S \cap M_e$ is characterized by one of the
following:

{\rm(1)}\qua $S \cap M_e$ is a collection of incompressible annuli (including
spanning annuli and possibly boundary parallel annuli) or is obtained
from such a collection by ambient 1--surgery along an arc which is
isotopic into $\partial M_e$.

{\rm(2)}\qua $M_e$ is homeomorphic to $(torus) \times I$ and there is a pair of
   simple closed curves $c, c' \subset (torus)$ such that $c \cap c'$
   consists of a single point $p \in (torus)$ and either $V \cap
   ((torus) \times I)$ or $W \cap ((torus) \times I)$ is a collar of
   $(c \times \{0\}) \cup (p \times I) \cup (c' \times \{1\})$.
\end{thm}

In the general case we can say the following:

\begin{thm} \label{thm:main}
Let $M$ be a totally orientable graph manifold.  Let $M = V \cup_S W$
be an irreducible Heegaard splitting.  Let $M = (V_1 \cup_{S_1} W_1)
\cup_{F_1} (V_2 \cup_{S_2} W_2) \cup_{F_2} \dots \cup_{F_{m-1}} (V_m
\cup_{S_m} W_m)$ be a weak reduction of $M = V \cup_S W$.  Set $M_i =
V_i \cup W_i$.  Then $M_i$ is a totally orientable generalized graph
manifold and $M_i = V_i \cup_{S_i} W_i$ is a strongly irreducible
Heegaard splitting.  In particular, $M_i = V_i \cup_{S_i} W_i$ is
standard.

Here $\chi(S) = \sum_i (\chi(S_i) - \chi(F_i))$.
\end{thm}

\begin{thm} \label{thm:postmain}
Let $M$ be a totally orientable graph manifold.  If $M = V \cup_S W$
is an irreducible Heegaard splitting, then it is the amalgamation of
standard Heegaard splittings of generalized graph submanifolds of $M$.
\end{thm}

The graph manifolds considered here are totally orientable, that is,
they are orientable $3$--manifolds and for each vertex manifold the
underlying surface of the orbit space is orientable.  It follows from
\cite[VI.34]{J} that an incompressible surface can be isotoped to be
either horizontal or vertical in each vertex manifold of a totally
orientable graph manifold.  In conjunction with the notion of
untelescoping a weakly reducible Heegaard splitting into a strongly
irreducible generalized Heegaard splitting, this observation reduces
the investigation at hand to the investigation of strongly irreducible
Heegaard splittings of generalized graph manifolds (for definitions,
see below).

In the investigation of strongly irreducible Heegaard splittings of
generalized graph manifolds, the nice properties of strongly
irreducible Heegaard splittings often reduce this investigation to a
study of the behaviour of the Heegaard splittings near the
characteristic submanifolds.  In this context, a theorem of D Cooper
and M Scharlemann completes the description of this behaviour, see
Proposition \ref{prop:acs} and Proposition \ref{prop:tcs}.  This
theorem may be found in \cite[Theorem 4.2]{CS}.

The theorem here is purely structural in the sense that it describes
the various ways in which a Heegaard splitting can be constructed.
Specifically, there are finitely many possible constructions.  Thus a
totally orientable graph manifold possesses only finitely many
Heegaard splittings up to isotopy (and hence also up to
homeomorphism).  It would be possible to extract a formula for the
genera of these Heegaard splittings.  However, this formula would be
long, cumbersome and not very enlightening.  But note that, in
particular, the program here enables a computation of Heegaard genus,
ie, the smallest possible genus of a Heegaard splitting, for totally
orientable graph manifolds.  To compute this genus, one need merely
consider the finitely many possible constructions, compute the
corresponding Euler characteristics, and find the extremal value.
This line of thought is pursued in \cite{SW}, where the genus of a
certain class of totally orientable graph manifolds is compared to the
rank, ie, the least number of generators, of the fundamental group
of these manifolds.

The theorem leaves open the question of classification.  There will be
some, though probably not too many, cases in which the various
constructions are isotopic.  More interestingly, there may be larger
scale isotopies.  Ie, there may be two Heegaard splittings of a
graph manifold that are isotopic but not via an isotopy fixing their
intersection with the decomposing tori.  Clearly, this leaves much
room for further investigation.

The global strategy is as follows: By Theorem \ref{thm:ind}, a
Heegaard splitting is the amalgamation of the strongly irreducible
Heegaard splittings arising in any of its weak reductions.  Thus, one
begins with a Heegaard splitting of a graph manifold.  One then
considers a weak reduction of this Heegaard splitting.  Cutting along
the incompressible surfaces in the weak reduction yields generalized
graph manifolds with strongly irreducible Heegaard splittings.  One
analyzes the possible strongly irreducible Heegaard splittings of
generalized graph manifolds.  Finally, one considers all possibilities
arising in the amalgamation of strongly irreducible Heegaard
splittings of generalized graph manifolds.

I wish to thank the many colleagues who have reminded me that a
complete report on this investigation is past due.  Among these are
Ian Agol, Hugh Howards, Yoav Moriah, Marty Scharlemann, Yo'av Rieck,
Eric Sedgwick and Richard Weidmann.  I also wish to thank the
MPIM-Bonn where part of this work was done.  This work was supported
in part by the grant NSF-DMS 0203680.

\section{Totally orientable graph manifolds}

For standard definitions pertaining to knot theory see for instance
\cite{BZ}, \cite{L} or \cite{R}.  For $3$--manifolds see \cite{He} or
\cite{J}.  Note that the terminology for graph manifolds has not been
standardized.

\begin{defn}
A Seifert manifold is a compact 3--manifold that admits a foliation by
circles.
\end{defn}

For a more concrete definition, see for instance \cite{J}.  The fact
that the simple definition here is in fact equivalent to more
concrete definitions follows from \cite{E}.  

\begin{defn}
The circles in the foliation of a Seifert fibered space $M$ are called
fibers.  The natural projection that sends each fiber to a point is
denoted by $p\co  M \rightarrow Q$.  The quotient space $Q$ is called the
base orbifold of $M$.  A fiber $f$ is called an exceptional fiber if
nearby fibers wind around $f$ more than once.  Otherwise, $f$ is
called a regular fiber.  The image under $p$ of a regular fiber is
called a regular point and the image under $p$ of an exceptional fiber
is called an exceptional point.
\end{defn}

The base orbifold is in fact a surface.  This follows from standard
facts about foliations in conjunction with \cite{E}.  It also follows
that there will be only finitely many exceptional fibers.

\begin{defn}
For $Y$ a submanifold of $X$, we denote an open regular neighborhood
of $Y$ in $X$ by $\eta(Y, X)$, or simply by $\eta(Y)$, if there is no
ambiguity concerning the ambient manifold.  Similarly, we denote a
closed regular neighborhood by $N(Y, X)$, or simply by $N(Y)$, if
there is no ambiguity concerning the ambient manifold.
\end{defn}

\begin{defn}
A surface $S$ in a Seifert fibered space $M$ is vertical if it
consists of fibers.  It is horizontal if it intersects all fibers
transversely.  It is pseudohorizontal if there is a fiber $f
\subset M$ such that $S \cap (M \backslash \eta(f))$ is horizontal and
$S \cap N(f)$ is a collar of $f$.
\end{defn}

It follows that a horizontal surface in a Seifert fibered space $M$
orbifold covers the base orbifold of $M$.

\begin{defn}
A Seifert fibered space is totally orientable if it is orientable as a
3--manifold and has an orientable base orbifold.
\end{defn}

We are now ready to define graph manifolds.

\begin{defn}
A graph manifold is a 3--manifold $M$ modelled on a finite graph
$\Gamma$ as follows:

{\rm(1)}\qua Each vertex $v$ of $\Gamma$ corresponds to a Seifert
fibered space, denoted by $M_v$ and called a vertex manifold;

{\rm(2)}\qua Each edge $e$ of $\Gamma$ corresponds to a 3--manifold homeomorphic
to $(torus) \times S^1$, denoted by $M_e$ and called an edge
manifold; 

{\rm(3)}\qua If an edge $e$ is incident to a vertex $v$, then this
incidence is realized by an identification of a boundary component of
$M_e$ with a boundary component of $M_v$ via a homeomorphism.

A graph manifold is totally orientable if each vertex manifold is
totally orientable.

The union of edge manifolds in $M$ is also called the characteristic
submanifold of $M$.  It is denoted by $\cal E$.  The image of a
boundary component of the characteristic submanifold of $M$ is a torus
called a decomposing torus.  It is denoted by $\cal T$.
\end{defn}

\begin{figure}[ht!]
\centerline{\epsfxsize=2.0in \epsfbox{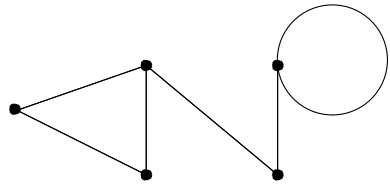}}
\caption{A model graph for a graph manifold}
\label{gm8.eps}
\end{figure}

A decomposing torus is, of course, also the image of a boundary
component of a vertex manifold.  But the converse is not always true.

\begin{rem}
We have placed no restrictions on the homeomorphism that identifies a
boundary component of an edge manifold with a boundary component of a
vertex manifold.  Thus according to this definition, there will be
Seifert fibered spaces that admit a description as a graph manifold
with non empty characteristic submanifold.  From the point of view of
the investigation here, this is often a useful way to think of such a
Seifert fibered space.  See \cite{W}.
\end{rem}

\begin{defn}
A boundary component of a vertex manifold $M_v$ of a graph manifold
$M$ that is also a boundary component of $M$ is called an exterior
boundary component of $M_v$.  We denote the union of exterior boundary
components of $M_v$ by $\partial_E M_v$.
\end{defn}

\section{Untelescoping and amalgamation}

We here give the basic definitions concerning Heegaard splittings,
strongly irreducible Heegaard splittings, untelescopings and
amalgamations.  Theorem \ref{thm:ind} below is crucial to the global
strategy employed in our investigation.

\begin{defn} \label{defn:cb} 
A \emph{compression body} is a $3$--manifold $W$ obtained from a
connected closed orientable surface $S$ by attaching $2$--handles to $S
\times \{0\} \subset S \times I$ and capping off any resulting
$2$--sphere boundary components.  We denote $S \times \{1\}$ by
$\partial_+W$ and $\partial W \backslash \partial_+W$ by
$\partial_-W$.  Dually, a compression body is a connected orientable
$3$--manifold obtained from a (not necessarily connected) closed
orientable surface $\partial_-W \times I$ or a $3$--ball by attaching
$1$--handles.

In the case where $\partial_-W = \emptyset$ (ie, in the case where a
$3$--ball was used in the dual construction of $W$), we also call $W$ a
\emph{handlebody}.  If $W = \partial_-W \times I$, we say that
$W$ is a \emph{trivial compression body}.
\end{defn}

\begin{defn}
A spine of a compression body $W$ is a 1--complex $X$ such that $W$
collapses to $\partial_-W \cup X$.
\end{defn}

\begin{defn} A \emph{set of defining disks} for a compression body
$W$ is a set of disks $\{D_1, \dots, D_n\}$ properly imbedded in $W$
with $\partial D_i \subset \partial_+W$ for $i = 1$, $\dots, n$ such
that the result of cutting $W$ along $D_1 \cup \dots \cup D_n$ is
homeomorphic to $\partial_-W \times I$ or to a $3$--ball in the case
that $W$ is a handlebody.
\end{defn}

\begin{defn} \label{defn:Heegaard splitting} A \emph{Heegaard
splitting} of a $3$--manifold $M$ is a decomposition $M = V \cup_S W$
in which $V$, $W$ are compression bodies, $V \cap W = \partial_+V =
\partial_+W = S$ and $M = V \cup W$.  We call $S$ the
\emph{splitting surface} or \emph{Heegaard surface}.
\end{defn}

The notion of strong irreducibility of a Heegaard splitting was
introduced by A. Casson and C. McA. Gordon in \cite{CG} and has proven
extremely useful.  

\begin{defn} \label{defn:wr} 
A Heegaard splitting $M = V \cup_S W$ is \emph{strongly
irreducible} if for any pair of essential disks $D \subset V$ and $E
\subset W$, $\partial D \cap \partial E \neq \emptyset$.
\end{defn}

Recall also the following related definitions:

\begin{defn}
A Heegaard splitting $M = V \cup_S W$ is \emph{reducible} if
there exists a pair of essential disks $D \subset V$ and $E \subset W$
such that $\partial D = \partial E$.  If $M = V \cup_S W$ is not
reducible, then it is \emph{irreducible}.

A Heegaard splitting $M = V \cup_S W$ is \emph{stabilized} if
there exists a pair of essential disks $D \subset V$ and $E \subset W$
such that $|\partial D \cap \partial E| = 1$.
\end{defn}

Though all compact $3$--manifolds admit Heegaard splittings, many do
not admit strongly irreducible Heegaard splitting.  This fact prompted
M. Scharlemann and A. Thompson to introduce the following notion of
generalized Heegaard splittings.

\begin{defn}
\label{general} A \emph{generalized Heegaard splitting} of a compact
orientable $3$--manifold $M$ is a decomposition $M = (V_1 \cup_{S_1}
W_1) \cup_{F_1} (V_2 \cup_{S_2} W_2) \cup_{F_2} \dots \cup_{F_{m-1}}
(V_m \cup_{S_m} W_m)$ such that each of the $V_i$ and $W_i$ is a union
of compression bodies with $\partial_+V_i = S_i = \partial_+W_i$ and
$\partial_-W_i = F_i = \partial_-V_{i+1}$.

We say that a generalized Heegaard splitting is \emph{strongly
irreducible} if each Heegaard splitting of a component of $M_i = V_i
\cup_{S_i} W_i$ is strongly irreducible and each $F_i$ is
incompressible in $M$.  We will denote $\cup_i F_i$ by $\cal F$ and
$\cup_i S_i$ by $\cal S$.  The surfaces in $\cal F$ are called the
\emph{thin levels} and the surfaces in $\cal S$ the
\emph{thick levels}.

Let $M = V \cup_S W$ be an irreducible Heegaard splitting.  We may
think of $M$ as being obtained from $\partial_-V \times I$ by
attaching all $1$--handles in $V$ (dual definition of compression body)
followed by all $2$--handles in $W$ (standard definition of compression
body), followed, perhaps, by $3$--handles.  An
\emph{untelescoping} of $M = V \cup_S W$ is a rearrangement of
the order in which the $1$--handles of $V$ and the $2$--handles of $W$
are attached yielding a generalized Heegaard splitting.  A
\emph{weak reduction} of $M = V \cup_S W$ is a strongly
irreducible untelescoping of $M = V \cup_S W$.
\end{defn}

Note that a weak reduction of a strongly irreducible Heegaard
splitting would just be the strongly irreducible Heegaard splitting
itself.  The Main Theorem in \cite{ST} implies the following:

\begin{thm} \label{thm:ST}
Let $M$ be an irreducible $3$--manifold.  Any Heegaard splitting $M = V
\cup_S W$ has a weak reduction.
\end{thm}

\begin{defn}
Let $N, L$ be $3$--manifolds with $R$ a closed subsurface of $\partial
N$, and $S$ a closed subsurface of $\partial L$, such that $R$ is
homeomorphic to $S$ via a homeomorphism $h$.  Further, let $(U_1,
U_2), (V_1, V_2)$ be Heegaard splittings of $N, L$ such that $N(R)
\subset U_1, N(S) \subset V_1$.  Then, for some $R' \subset
\partial N \backslash R$ and $S' \subset \partial L \backslash S$,
$U_1 = N(R \cup R') \cup (1-handles)$ and $V_1 = N(S \cup S')
\cup (1-handles)$.  Here $N(R)$ is homeomorphic to $R \times I$ via
a homeomorphism $f$ and $N(S)$ is homeomorphic to $S \times I$ via
a homeomorphism $g$.  Let $\sim$ be the equivalence relation on $N
\cup L$ generated by

(1)\qua $x \sim y$ if $x, y \epsilon \eta(R)$ and $p_1 \cdot f(x) = p_1
    \cdot f(y)$,

(2)\qua  $x \sim y$ if $x, y \epsilon \eta(S)$ and $p_1 \cdot g(x) = p_1
    \cdot g(y)$,

(3)\qua  $x \sim y$ if $x \epsilon R$, $y \epsilon S$ and $h(x) = y$,

where $p_1$ is projection onto the first coordinate.  Perform
isotopies so that for $D$ an attaching disk for a $1$--handle in $U_1,
D'$ an attaching disk for a $1$--handle in $V_1$, $[D] \cap [D'] =
\emptyset$.  Set $M = (N \cup L)/\sim, W_1 = (U_1 \cup V_2)/\sim,$ and
$W_2 = (U_2 \cup V_1)/\sim$.  In particular, $(N(R) \cup N(S)/\sim)
\cong R, S$.  Then $W_1 = V_2 \cup N(R') \cup (1--handles)$, where the
$1$--handles are attached to $\partial_+V_2$ and connect $\partial
N(R')$ to $\partial_+V_2$, and hence $W_1$ is a compression body.
Analogously, $W_2$ is a compression body.  So $(W_1, W_2)$ is a
Heegaard splitting of $M$.  The splitting $(W_1, W_2)$ is called the
\emph{amalgamation of $(U_1, U_2)$ and $(V_1, V_2)$ along {R, S}
via $h$}.
\end{defn}

Theorem \ref{thm:ST} together with \cite[Proposition 2.8]{Sch} implies
the following:

\begin{thm} \label{thm:ind} 
Suppose $M = V \cup_S W$ is an irreducible Heegaard splitting and $M =
(V_1 \cup_{S_1} W_1) \cup_{F_1} (V_2 \cup_{S_2} W_2) \cup_{F_2} \dots
\cup_{F_{m-1}} (V_m \cup_{S_m} W_m)$ a weak reduction of $M = V \cup_S
W$.  Then the amalgamation of $M = (V_1 \cup_{S_1} W_1) \cup_{F_1}
(V_2 \cup_{S_2} W_2) \cup_{F_2} \dots \cup_{F_{m-1}} (V_m \cup_{S_m}
W_m)$ along $F_1 \cup \dots \cup F_{m-1}$ is $M = V
\cup_S W$.
\end{thm}

One of the nice properties of strongly irreducible Heegaard splittings
is apparent in the following lemma which is a deep fact and is proven,
for instance, in \cite[Lemma 6]{Sch4}.

\begin{lem}  \label{lem:essential} 
Suppose $M = V \cup_S W$ is a strongly irreducible Heegaard splitting
and $P \subset M$ an essential incompressible surface.  Then $S$ can
be isotoped so that $S \cap P$ consists only of curves essential in
both $S$ and $P$.
\end{lem}

\section{Incompressible surfaces and generalized graph\nl manifolds}

In the arguments that follow, we employ the ideas of untelescoping and
amalgamation.  In this section, we describe the incompressible
surfaces that arise in a weak reduction of a Heegaard splitting.  We
then describe the 3--manifolds that result from cutting a graph
manifold along such incompressible surfaces.  We will call these
3--manifolds generalized graph manifolds.  Later, we will consider the
strongly irreducible Heegaard splittings on these generalized graph
manifolds.  

\begin{rem} \label{inc}
An edge manifold of a totally orientable graph manifold $M$ is
homeomorphic to $(torus) \times I$.  There are infinitely many
distinct foliations of $(torus) \times I$ as an annulus bundle over
the circle.  The incompressible surfaces in $(torus) \times I$ are
tori isotopic to $(torus) \times \{point\}$, annuli isotopic to the
annular fibers in the foliations of $(torus) \times I$ as an annulus
bundle over the circle and annuli parallel into $\partial ((torus)
\times I)$.
\end{rem}

\begin{lem} \label{lem:inc}
Let $F$ be an incompressible surface in a totally orientable graph
manifold $M$.  Then $F$ may be isotoped so that in each edge manifold
it consists of incompressible tori and essential annuli and in each
vertex manifold it is either horizontal or vertical.
\end{lem}

\proof 
Let $\cal T$ be the collection of decomposing tori for $M$.
Since $F$ and $\cal T$ are incompressible, $F$ may be isotoped so that
$F \cap \cal T$ consists only of curves essential in both $F$ and
$\cal T$.  We may assume that this has been done in such a way that
the number of components in $F \cap \cal T$ is minimal.  Let $N$ be a
component of $M \backslash \cal T$, then $F \cap N$ is incompressible.
Furthermore, no component of $F \cap N$ is an annulus parallel into
$\cal T$.

Suppose $F \cap N$ is boundary compressible in $N$.  Let $\hat D$ be a
boundary compressing disk for $F \cap N$.  Then $\partial \hat D = a
\cup b$, with $a \subset \partial N$ and $b \subset F$.  Since $F \cap
\cal T$ consists only of curves essential in both $F$ and $\cal T$,
the component $A$ of $\partial N \backslash (F \cap \partial N)$ that
contains $a$ is an annulus.  Let $B(\hat D)$ be a bicollar of $\hat
D$.  Then $\partial B(\hat D)$ has two components, $\hat D_0, \hat
D_1$.  Consider the disk $D = (A \backslash (A \cap B(\hat D)) \cup
\hat D_0 \cup \hat D_1$.  Since $F$ is incompressible, $D$ must be
parallel to a disk in $F$, but this implies that the number of
components of $F \cap \cal T$ is not minimal, a contradiction.  Thus,
$F \cap N$ is boundary incompressible in $N$.

If $N$ is an edge manifold, then $F \cap N$ is as required by Remark
\ref{inc}.  If $N$ is a vertex manifold, then \cite[VI.34]{J} allows
three possibilities for $F \cap N$: (1) $F \cap N$ is vertical; (2) $F
\cap N$ is horizontal; or (3) $F \cap N$ is the boundary of a twisted
$I$--bundle over a horizontal surface $\hat F^N$ in $N$.  For a
boundary incompressible surface in $N$ this latter possibility would
imply that there is a nonorientable horizontal surface $\hat F^N$ in
$C$.  In particular, $\hat F^N$ would be a cover of the base orbifold.
But this is impossible.  Hence $F \cap N$ is either horizontal or
vertical.  Hence $F \cap N$ is as required.  \endproof

The following definition describes the 3--manifolds that result when a
totally orientable graph manifold is cut along incompressible
surfaces.

\begin{defn}
A generalized graph manifold is a 3--manifold $M$ modelled on a finite
graph $\Gamma$ as follows:

{\rm(1)}\qua Each vertex $v$ of $\Gamma$ corresponds either to a Seifert fibered
space or to a 3--manifold homeomorphic to $(compact \; \surface) \times
[0, 1]$.  This manifold is denoted by $M_v$ and called a vertex
manifold.

{\rm(2)}\qua Each edge $e$ of $\Gamma$ corresponds either to a 3--manifold
homeomorphic to $(torus) \times I$ or to a 3--manifold homeomorphic
to $(annulus) \times I$.  This manifold is denoted by $M_e$ and
called an edge manifold.

{\rm(3)}\qua If the edge manifold $M_e$ is homeomorphic to $(torus) \times I$
and $e$ is incident to a vertex $v$, then this incidence is realized
by an identification of a boundary component of $M_e$ with a boundary
component of $M_v$.  In particular, $M_v$ must be Seifert fibered.

{\rm(4)}\qua If the edge manifold $M_e$ is homeomorphic to $(annulus) \times I$
and $e$ is incident to a vertex $v$, then this incidence is realized
by an identification of a component of $(\partial (annulus)) \times
I$, with a subannulus of $\partial M_v$.  If $M_v$ is Seifert fibered,
then this subannulus of $\partial M_v$ consists of fibers of $M_v$.
If $M_v$ is homeomorphic to $(compact \; \surface) \times I$, then this
subannulus is a component of $(\partial (compact \; \surface)) \times
I$.

{\rm(5)}\qua If a vertex manifold $M_v$ is homeomorphic to $(compact \; \surface)
   \times I$, then the valence of $v$ equals the number of components
   of $(\partial (compact \; \surface)) \times I$.  Ie, each
   component of $(\partial (compact \; \surface)) \times I$ is
   identified with a subannulus of the boundary of an edge manifold.

A generalized graph manifold is totally orientable if
each vertex manifold that is Seifert fibered is totally orientable and
each vertex manifold that is not Seifert fibered is homeomorphic to
$(compact \; orientable \; \surface) \times I$.

The union of edge manifolds in $M$ is also called the characteristic
submanifold of $M$.  It is denoted by $\cal E$.  The image of a torus
or annulus, respectively, along which an indentification took place is
called a decomposing torus or decomposing annulus, respectively.  The
union of decomposing tori and annuli is denoted by $\cal T$.
\end{defn}

Consider the case in which an edge manifold $M_e = (torus) \times
[0,1]$ of a graph manifold is cut along an incompressible torus $T =
(torus) \times \{point\}$.  When $M_e$ is cut along $T$, the remnants
of $M_e$ are $(torus) \times [0, \frac{1}{2}]$ and $(torus) \times
[\frac{1}{2}, 1]$.  Each of these remnants forms a collar of a vertex
manifold.  A foliation of $(torus)$ by circles may be chosen in such a
way that the Seifert fibration of the vertex manifold extends across
the remnant.  We may thus ignore these remnants, as we do in the above
definition of generalized graph manifolds.  This facilitates the
discussion of strongly irreducible Heegaard splittings of generalized
graph manifolds.  Later, when considering amalgamations of strongly
irreducible Heegaard splittings of generalized graph manifolds, we
will have to reconsider these remnants.

\begin{defn}
A portion of the boundary of a vertex manifold $M_v$ of a generalized
graph manifold $M$ that is contained in $\partial M$ is called an
exterior boundary component of $M_v$.  We denote the union of exterior
boundary components of $M_v$ by $\partial_E M_v$.
\end{defn}

\section{Motivational examples of Heegaard splittings}

In this section we describe some examples of Heegaard splittings for
graph manifolds and generalized graph manifolds.  The following
definition facilitates describing the structure of certain surfaces.  

\begin{defn}
Let $F$ be a surface in a $3$--manifold $M$ and $\alpha$ an arc with
interior in $M \backslash F$ and endpoints on $F$.  Let $C(\alpha)$ be
a collar of $\alpha$ in $M$.  The boundary of $C(\alpha)$ consists of
an annulus $A$ together with two disks $D_1, D_2$, which we may assume
to lie in $F$.  We call the process of replacing $F$ by $(F \backslash
(D_1 \cup D_2)) \cup A$ \emph{performing ambient $1$--surgery on
$F$ along $\alpha$}.
\end{defn}

The process of ambient $1$--surgery on a surface along an arc is
sometimes informally referred to as ``attaching a tube''.

\begin{ex}
Let $Q$ be a closed orientable surface.  The standard Heegaard
splitting of $Q \times \mathbb S^1$ may be constructed in more than
one way.  In particular, consider a small disk $D \subset Q$ and a
collection $\Gamma$ of arcs that cut $Q \backslash D$ into a disk.
Let $S$ be the result of performing ambient $1$--surgery on $\partial D
\times \mathbb S^1$ along $\Gamma \times \{point\}$.  Then $S$ is the
splitting surface of a Heegaard splitting $Q \times \mathbb S^1 = V
\cup_S W$ (for details, \cite{Sch}).  See Figure \ref{gm9.ps}.  One of
the handlebodies is $((shaded \; disk) \times \mathbb S^1) \cup
N(dashed \; arcs)$.

\begin{figure}[ht!]
\centerline{\epsfxsize=2.0in \epsfbox{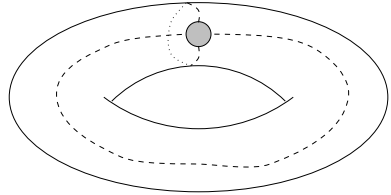}}
\caption{Schematic for Heegaard splitting of $Q \times \mathbb S^1$
\sl}
\label{gm9.ps}
\end{figure}

The same Heegaard splitting of $Q \times \mathbb S^1$ may be obtained
in another way: Partition $\mathbb S^1$ into two intervals $I_1, I_2$
that meet in their endpoints.  Consider two distinct points $p, q \in
Q$.  Then the surface obtained by performing ambient 1--surgery on $Q
\times (I_1 \cap I_2)$ along $(p \times I_1) \cup (q \times I_2)$ is
isotopic to $S$ above.
\end{ex}

The fact that this Heegaard splitting can be constructed either from a
vertical torus or from horizontal surfaces via ambient 1--surgery is
likely to be a very special feature.  But there are no general
techniques for detecting the sort of global isotopies that allow this
to happen.

Note also that this Heegaard splitting is not strongly irreducible.
The two descriptions of this Heegaard splitting hint at distinct weak
reductions.  In one of these weak reductions, the incompressible
surfaces would consist of vertical tori.  In the other, the
incompressible surfaces would consist of horizontal incompressible
surfaces.  In weak reductions of the latter type, we see that the
Heegaard splitting is an amalgamation of two Heegaard splittings with
mostly horizontal splitting surface.

The Heegaard splitting described turns out to be the only irreducible
Heegaard splitting for a manifold of the form $(closed \; orientable
\; \surface) \times \mathbb S^1$.  This fact is the main theorem of
\cite{Sch}.  The first description of the construction can be
generalized to Seifert fibered spaces to provide the canonical
Heegaard splittings for totally orientable Seifert fibered spaces, see
\cite{BO} and \cite{MS}.  The Heegaard splittings arising from this
construction have been termed vertical.  This terminology has created
some confusion, because the splitting surface of a vertical Heegaard
splitting is not vertical as a surface.  Here we will continually
focus on the splitting surface.  In particular, we will want to
distinguish between surfaces that are vertical and surfaces that are
the splitting surface of a vertical Heegaard splitting.  For this
reason, we will augment the existing terminology and refer to the
splitting surface of a vertical Heegaard splitting as
``pseudovertical''.

We recall the definition of a vertical Heegaard splitting for a
Seifert fibered space.  The structure of Heegaard splittings for
totally orientable Seifert fibered spaces has been completely
described in \cite{MS}.  Thus we may restrict our attention to
Seifert fibered spaces with non empty boundary in the definition
below.  

\begin{figure}[ht!]
\centerline{\epsfxsize=2.0in \epsfbox{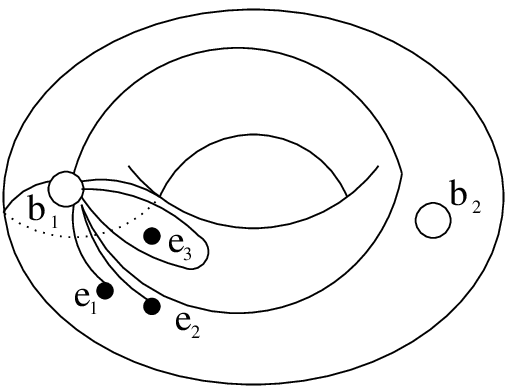}}
\caption{Schematic for a vertical Heegaard splitting of a Seifert
fibered space \sl}
\label{gm12a.ps}
\end{figure}

\begin{defn}
Let $M$ be a Seifert fibered space with $\partial M \neq
\emptyset$. Denote the base orbifold of $M$ by $O$.  Denote the
exceptional fibers of $M$ by $f_1, \dots, f_n$ and the corresponding
exceptional points in $O$ by $e_1, \dots, e_n$.  Denote the boundary
components of $M$ by $B_1, \dots, B_m$ and the corresponding boundary
components of $O$ by $b_1, \dots, b_m$,

Partition $f_1, \dots, f_n$ into two subsets: $f_1, \dots, f_i$, the
fibers that will lie in $V$ and $f_{i+1}, \dots, f_n$, the fibers that
will lie in $W$.  Then partition $B_1, \dots, B_m$ into two subsets:
$B_1, \dots, B_j$, the boundary components that will lie in $V$ and
$B_{j+1}, \dots, B_m$, the boundary components that will lie in $W$.

We may assume that $j \geq 1$, for otherwise we may interchange the
roles of $V$ and $W$ in the construction below.  Let $\Gamma$ be a
collection of arcs in $O$ each with at least one endpoint on $b_1$
such that $O \backslash \Gamma$ is a regular neighborhood of $e_{i+1}
\cup \dots \cup e_n \cup b_{j+1} \cup \dots \cup b_m$ or of a point,
if this set is empty.  See Figure \ref{gm12a.ps}.

Set $V = N(f_1 \cup \dots \cup f_i \cup B_1 \cup \dots \cup B_j \cup
\Gamma)$ and set $W = closure(M \backslash V)$.  Set $S = \partial_+V
= \partial_+W$.  Then $M = V \cup_S W$ is a Heegaard splitting.  (For
details, see \cite{Sch1} or \cite{MS}.)  A Heegaard splitting of a
Seifert fibered space with non empty boundary constructed in this
manner is called a vertical Heegaard splitting.  A pseudovertical
surface is a surface that is the splitting surface of a vertical
Heegaard splitting.
\end{defn}

A little more work is required to extend this notion to the setting of
graph manifolds.  Here we consider the intersection of a Heegaard
splitting of a generalized graph manifold with a vertex manifold that
is a Seifert fibered space.  Denote the Heegaard splitting by $M = V
\cup_S W$ and the vertex manifold by $M_v$.  Here $S \cap M_v$ is not
necessarily connected.  Furthemore, $S \cap M_v$ has boundary.  A
concrete description of all possible such surfaces would be quite
extensive.  For this reason, we give the following, purely structural,
definition:

\begin{defn}
Let $M = V \cup_S W$ be a Heegaard splitting of a generalized graph
manifold with non empty characteristic submanifold.  Let $M_v$ be a
Seifert fibered vertex manifold of $M$.  We say that $S \cap M_v$ is
pseudovertical if the following holds:

There is a collection of vertical annuli and tori ${\cal A} \subset
M_v$. And there is a collection of arcs $\Gamma$ in the interior of
$M_v$ such that each endpoint of each arc in $\Gamma$ lies in ${\cal
A}$ and such that $\Gamma$ projects to a collection of disjoint
imbedded arcs in $O_v$.  And $S \cap M_v$ is obtained from ${\cal A}$
by ambient 1--surgery along $\Gamma$.
\end{defn}

The assumption that $M = V \cup_S W$ is a Heegaard splitting places
strong restrictions on $S$.  If $S \cap M_v$ is pseudovertical, then
$S \cap (\partial M_v \backslash \partial_E M_v)$ consists of vertical
curves.  Thus $V \cap (\partial M_v \backslash \partial_E M_v)$ and $W
\cap (\partial M_v \backslash \partial_E M_v)$ consist of annuli.
Each such annulus is either a spanning annulus in $V$ or $W$,
or it has both boundary components in $S$.

Consider the result of cutting a compression body $V$ along an annulus
$A$ with $\partial A \subset \partial_+V$.  If the annulus is
inessential, then the effect is nil.  If the annulus is essential,
then the result is again a, possibly disconnected, compression body,
see \cite[Lemma 2]{Sch4}.  If the annulus is a spanning annulus, then
the result is a handlebody.  But in this context we should think of it
as a compact 3--manifold of the form $((compact\; \surface) \times I)
\cup (1-handles)$.  The only way this can happen, given the structure
of $S \cap M_v$, is if the components of $V \cap M_v$ and $W \cap M_v$
are constructed from vertical solid tori and perhaps components
homeomorphic to $(annulus) \times \mathbb S^1$ by attaching
``horizontal'' 1--handles.  For more concrete computations, see
\cite{SW}.

It is a non trivial fact that for Seifert fibered spaces the two
definitions of pseudovertical surfaces coincide.  This follows from
\cite[Proposition 2.10]{Sch} via Lemma \ref{lem:bcore}
(an adaptation of the central argument in \cite{Sch1}) along with
Lemma \ref{lem:ex}.  For an illustration, see the final remarks in the
example below.

\begin{ex} \label{ex:c}
Let $M$ be a Seifert fibered space with base orbifold a disk and with
two exceptional fibers $f_1, f_2$.  Let $T$ be a boundary parallel
torus and let $\alpha$ be an arc connecting $T$ to itself that
projects to an imbedded arc that separates the two exceptional points.
Let $S$ be the result of performing ambient $1$--surgery on $T$ along
$\alpha$.  Then $S$ is the splitting surface of a Heegaard splitting
of $M = V \cup_S W$ (see \cite{MS}).

Now consider two copies $M_1, M_2$ of $M$ with Heegaard splittings
$M_i = V_i \cup_{S_i} W_i$.  We may assume that $\partial M_1 \subset
V_1$ and $\partial M_2 \subset W_2$.  We identify $\partial M_1$ and
$\partial M_2$ to obtain a 3--manifold $\tilde M$.  Juxtaposing the two
Heegaard splittings provides a generalized strongly irreducible
Heegaard splitting.  It is indicated schematically in Figure
\ref{gm10a.ps}.  Amalgamating the two Heegaard splittings along
$\partial M_1, \partial M_2$ results in a Heegaard splitting $\tilde M
= V \cup_S W$.  The result is schematically indicated in Figure
\ref{gm11.ps}.  The circles correspond to vertical tori.  The
splitting surfaces of the Heegaard splittings are obtained by
performing ambient 1--surgery on these tori along arcs in $\tilde M$
corresponding to the dashed arcs.

\begin{figure}[ht!]
\centerline{\epsfxsize=1.8in \epsfbox{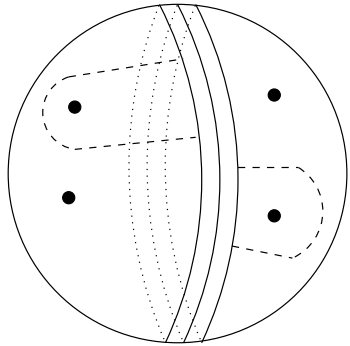}}
\caption{Schematic for a generalized strongly irreducible Heegaard splitting}
\label{gm10a.ps}
\end{figure}

\begin{figure}[ht!]
\centerline{\epsfxsize=1.8in \epsfbox{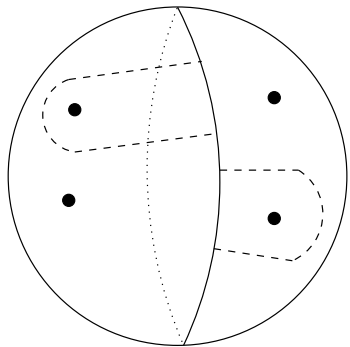}}
\caption{Schematic for Heegaard splitting after amalgamation}
\label{gm11.ps}
\end{figure}

The manifold $\tilde M$ obtained when $\partial M_1$ and $\partial
M_2$ are identified is a graph manifold modelled on a graph with two
vertices and one edge connecting the two vertices.  The vertex
manifolds are slightly shrunken versions of $M_1$ and $M_2$.  The edge
manifold is a collar of the image of $\partial M_1$ and $\partial M_2$
in $\tilde M$.

If the homeomorphism that identifies $\partial M_1$ with $\partial
M_2$ is fiberpreserving, then $\tilde M$ is in fact a Seifert fibered
space.  In this case $S$ is isotopic to the surface indicated
schematically in Figure \ref{gm11a.ps}.  This surface is obtained by
performing ambient 1--surgery along two arcs corresponding to the
dashed arcs on the two vertical tori corresponding to the two solid
circles.

\begin{figure}[ht!]
\centerline{\epsfxsize=1.8in \epsfbox{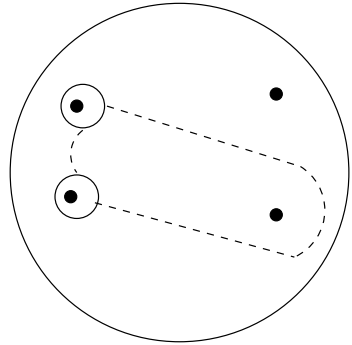}}
\caption{Schematic for Heegaard splitting of a Seifert fibered space}
\label{gm11a.ps}
\end{figure}

\end{ex}

The above construction can be generalized to arbitrary graph
manifolds.  The resulting Heegaard splittings can be considered the
canonical Heegaard splittings and generalized Heegaard splittings of
graph manifolds.  Contemplation of Figure \ref{gm11.ps} might lead
one to believe that amalgamations of strongly irreducible Heegaard
splittings still have enough structure to be described in the
terminology used here.  But this is not the case, as we shall see in
the following example.

\begin{ex} \label{ex:nostruc}
Let $M_1$ be as above.  Let $P$ be a thrice punctured $\mathbb
S^2$. Denote the boundary components of $P$ by $b_1, b_2, b_3$.  Let
$M_2$ be the 3--manifold obtained from $P \times \mathbb S^1$ by
identifying $b_2 \times \mathbb S^1$ and $b_3 \times \mathbb S^1$ via
a homeomorphism that takes $b_2 \times \{point\}$ to $\{point\} \times
\mathbb S^1$.  

Here $M_2$ is a graph manifold modelled on a graph with one vertex and
one edge.  It has a Heegaard splitting depicted schematically in
Figure \ref{gmns1.ps}.  In $P \times \mathbb S^1$, we take the
products of the regions pictured.  This yields a white $(annulus)
\times \mathbb S^1$ and a shaded solid torus.  When $b_2 \times
\mathbb S^1$ and $b_3 \times \mathbb S^1$ are identified, both
$(annulus) \times \mathbb S^1$ and the solid torus meet themselves in
(square) disks.  Thus we obtain a (strongly irreducible) Heegaard
splitting of genus 2 for $M_2$.  Denote the compression body by $V_2$
and the handlebody by $W_2$.

\begin{figure}[ht!]
\centerline{\epsfxsize=2.0in \epsfbox{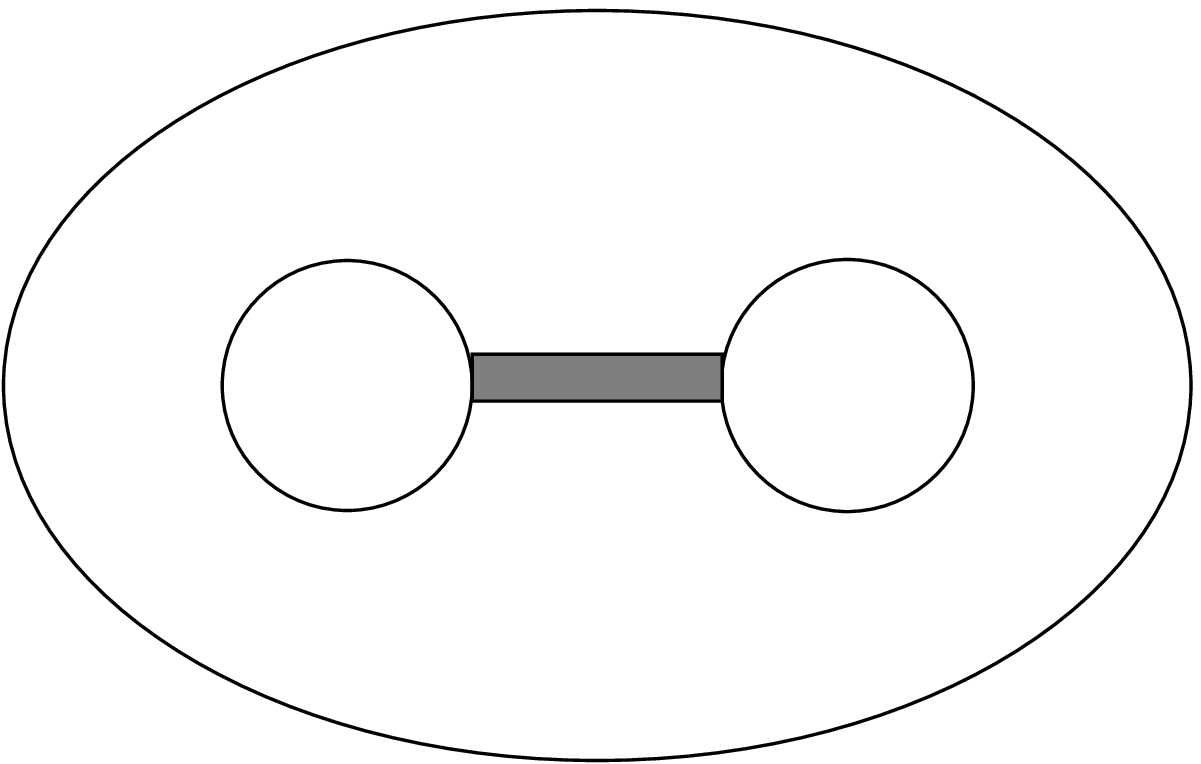}}
\caption{Schematic for Heegaard splitting of a graph manifold}
\label{gmns1.ps}
\end{figure}

This Heegaard splitting can also be constructed by ambient 1--surgery
on a boundary parallel torus along an arc as pictured in Figure
\ref{gmns2.ps}.

\begin{figure}[ht!]
\centerline{\epsfxsize=2.0in \epsfbox{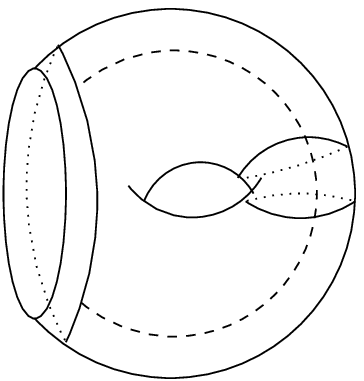}}
\caption{Schematic for Heegaard splitting of a graph manifold}
\label{gmns2.ps}
\end{figure}

The description indicated schematically in Figure \ref{gmns1.ps} is
preferable to the one indicated schematically in Figure
\ref{gmns2.ps}.  This is because the splitting surface of the former
is vertical in the Seifert fibered vertex manifold and has a very
special structure in the edge manifold.  The latter intersects the
edge manifold in a tube.  

Now let $M$ be a 3--manifold obtained by identifying the boundary
component of $M_1$ to the boundary component of $M_2$.  Then $M$ is a
graph manifold modelled on a graph as in Figure \ref{gmns3.ps}.

\begin{figure}[ht!]
\centerline{\epsfxsize=2.0in \epsfbox{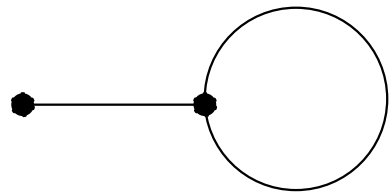}}
\caption{The graph on which $M$ is modelled}
\label{gmns3.ps}
\end{figure}

Here $M$ inherits a strongly irreducible generalized Heegaard
splitting that is the juxtaposition of the strongly irreducible
Heegaard splittings of $M_1$ and $M_2$.  If we wish to amalgamate the
two Heegaard splittings, recall that the collar neighborhood of
$\partial M_2$ that lies in $V_2$ will be identified to a single
torus.  To prevent this identification from interfering with the
decomposing tori of $M_2$, we are forced to consider the second
description of the Heegaard splitting $M_2 = V_2 \cup_{S_2} W_2$.  
But this means that our description of the resulting Heegaard 
splitting involves a tube which runs through an edge manifold.
\end{ex}

The tube in the above description of the Heegaard splitting does not
fit nicely into the characterization of Heegaard splittings as
horizontal, pseudohorizontal, vertical or pseudovertical.  The moral
is, that such tubes arise and can be quite complicated.  A more
thorough investigation of the nature of such tubes will be the subject
of a further investigation.  However, such tubes do not arise in
strongly irreducible Heegaard splittings.  For this reason we prefer
to think of our Heegaard splittings as amalgamations of strongly
irreducible Heegaard splittings.

The following two examples illustrate more peculiar Heegaard
splittings that arise under special circumstances.

\begin{ex}
Let $N$ be a Seifert fibered space with base orbifold the sphere and
with four exceptional fibers $f_1, \dots, f_4$ with carefully chosen
invariants: $\frac{1}{2}, \frac{1}{2}, \frac{1}{2}, \frac{l}{2l+1}$.
Here $N \backslash \eta(f_4)$ is a Seifert fibered manifold with
boundary and hence fibers over the circle.  More specifically, it
fibers as a once punctured torus bundle over the circle.  By
partitioning the circle into two intervals $I_1, I_2$ that meet in
their endpoints, we obtain a decomposition $N \backslash \eta(f_4) =
V' \cup_{S'} W'$ with $V' = (once\;punctured\;torus) \times I_1$, $W'
= (once\;punctured\;torus) \times I_2$ and $S' = two \;
once\;punctured\;tori$.

Note that this decomposition is not a Heegaard splitting in the sense
used here as $S'$ is not closed.  Now the carefully chosen invariants
guarantee that the boundary of a meridian disk of $N(f_4)$ meets
$\partial V'$ in a single arc.  In particular, $V = V' \cup N(f_4)$ is
also a handlebody (of genus two).  Setting $W = W'$, this defines a
Heegaard splitting $N = V \cup_S W$.  The splitting surface of this
Heegaard splitting is a horizontal surface away from $N(f_4)$.  And
after a small isotopy, $S \cap N(f_4)$ is a collar of $f_4$.  Thus $S$
is pseudohorizontal.
\end{ex}

\begin{ex}
Let $Q$ be a once punctured torus.  Set $M_i = Q \times \mathbb S^1$
for $i = 1, 2$.  Partition $\mathbb S^1$ into two intervals $I_1, I_2$
meeting in their endpoints.  Set $V_i = Q \times I_1 \subset M_i$ and
$W_i = Q \times I_2 \subset M_i$.  Let $T_i = \partial M_i$ and $c_i =
\partial Q \times \{point\} \subset T_i$.  Identify $T_1$ and $T_2$
via a homeomorphism so that $| c_1 \cap c_2 | = 1$.  Then $V_1$ and
$V_2$ meet in a (square) disk, hence $V = V_1 \cup V_2$ is a
handlebody.  Similarly, $W = W_1 \cup W_2$ is a handlebody.  Set $S =
\partial V = \partial W$, then $M = V \cup_S W$ is a Heegaard
splitting of $M = M_1 \cup_{T_1 = T_2} M_2$.

Here $M$ is a graph manifold modelled on a graph with two vertices and
one edge connecting the two vertices.  Its characteristic submanifold
is a collar of $T_1 = T_2$.  In the vertex manifolds, $S$ is
horizontal.  In the edge manifold, $S$ has a very specific structure.
\end{ex}

\begin{defn}
Let $M$ be a generalized graph manifold with characteristic
submanifold $\cal E$.  Let $M = V \cup_S W$ be a Heegaard splitting.
We say that $M = V \cup_S W$ is \emph{standard} if $S$ can be
isotoped so that for each vertex manifold $M_v$ of $M$, $S \cap M_v$
is either horizontal, pseudohorizontal, vertical or pseudovertical and
such that for each edge manifold $M_e$ of $M$, $S \cap M_e$ is
characterized by one of the following:

{\rm(1)}\qua $S \cap M_e$ is a collection of incompressible annuli (including
spanning annuli and possibly boundary parallel annuli) or is obtained
from such a collection by ambient 1--surgery along an arc which is
isotopic into $\partial M_e$.

{\rm(2)}\qua $M_e$ is homeomorphic to $(torus) \times I$ and there is a pair of
   simple closed curves $c, c' \subset (torus)$ such that $c \cap c'$
   consists of a single point $p \in (torus)$ and either $V \cap
   ((torus) \times I)$ or $W \cap ((torus) \times I)$ is a collar of
   $(c \times \{0\}) \cup (p \times I) \cup (c' \times \{1\})$.
\end{defn}

\section{The active component}

In this section we consider a generalized graph manifold $W$. We show
that if $M = V \cup_S W$ is a strongly irreducible Heegaard splitting,
then $S$ may be isotoped so that it is incompressible away from a
single vertex or edge manifold of $M$.

\begin{lem} \label{lem:C}
Let $M$ be a generalized graph manifold with characteristic
submanifold $\cal T$.  Let $M = V \cup_S W$ be a strongly irreducible
Heegaard splitting. Let ${\cal D}_V$ be a collection of defining disks
for $V$ and ${\cal D}_W$ a collection of defining disks for $W$.
There is a vertex or edge manifold $N$ of $M$ so that, after isotopy,
each outermost disk component of both ${\cal D}_V \backslash ({\cal T
\cap D}_V)$ and of ${\cal D}_W \backslash ({\cal T \cap D}_W)$ lies in
$N$.  Moreover, for each vertex or edge manifold $\tilde N \neq N$, $S
\cap \tilde N$ is incompressible.
\end{lem}

\proof
Isotope $S$ so that $\cal T$$\cap S$ consists only of curves
essential in both $\cal T$ and $S$.  Furthermore, assume that the
isotopy has been chosen so that $|\cal T$$\cap S|$ is minimal subject
to this condition.  Suppose that $D'$ is an outermost subdisk of
${\cal D}_V \backslash (\cal T \cap D)$.  Let $N$ be the vertex or
edge manifold of $M$ containing $D'$.  Then either $D'$ is a disk in
the interior of $V$ or $D'$ meets an annular component $A$ of $V \cap
\partial N$.  In case of the latter, $\partial D'$ meets $A$ in a
single arc, $a$.  Let $D''$ be the disk obtained by cutting $A$ along
$a$, adding two copies of $D' $ and isotoping the result to be a
properly imbedded disk in $V$.  The assumption that $|\cal T$$\cap S|$
be minimal guarantees that $D'' $ is an essential disk in $V$.  Thus
in both cases, there is an essential disk properly imbedded in $V$
that lies in the interior of $N$, we refer to this disk as $D$.

Similarly, consider an outermost subdisk of ${\cal D}_W \backslash
({\cal T} \cap {\cal D}_W)$.  The above argument shows that in a
vertex or edge manifold $N'$ of $M$, there is an essential disk $E$
properly imbedded in $W$.  Since $M = V \cup_S W$ is strongly
irreducible, $D$ must meet $E$, hence the vertex or edge manifold $N'$
must coincide with the vertex or edge manifold $N$.

It follows that for each vertex or edge manifold $\tilde N \neq N$, $S
\cap \tilde N$ is incompressible.
\endproof

\begin{defn}
If $M$ is a generalized graph manifold, with a strongly irreducible
Heegaard splitting $M = V \cup_S W$, then the vertex or edge manifold
$N$ as in Lemma \ref{lem:C} is called the \emph{active component
of $M = V \cup_S W$}.
\end{defn}

\section{What happens in the active component?}

The possibilities for the active component depend on the type of the
active component.  There are five possibilities.  We discuss each in
turn.

\subsection{Seifert fibered vertex manifold with exterior boundary}

We first consider the case in which the active component of $M = V
\cup_S W$ is a vertex manifold $M_v$ that is a Seifert fibered space
and that has exterior boundary.  This situation has been studied
extensively in the more restricted case in which $M$ itself is a
Seifert fibered space.

\begin{lem} \label{lem:bcore}
Suppose $M$ is a totally orientable Seifert fibered space with non
empty boundary and $M = V \cup_S W$ is a Heegaard splitting.  Then
each exceptional fiber of $M_v$ is a core of either $V$ or $W$.
\end{lem}

\proof
This is \cite[Lemma 4.1]{Sch1}, the central argument in \cite{Sch1}.
\endproof

\begin{thm} \label{thm:sch1}
Suppose $M$ is a totally orientable Seifert fibered space with non
empty boundary.  Suppose $M = V \cup_S W$ is a Heegaard splitting.
Then $S$ is pseudovertical.
\end{thm}

\proof
This is the main theorem of \cite{Sch1}.
\endproof

The argument extends to a more general setting.  Indeed, the isotopy
performed in \cite[Lemma 4.1]{Sch1} takes place within a small regular
neighborhood of a saturated annulus.  In particular, the isotopy can
here be performed entirely within $M_v$.  The only requirements on
this saturated annulus are that one boundary component has to lie on
the splitting surface of the Heegaard splitting and the other has to
wrap around the exceptional fiber.

\begin{lem} \label{lem:core}
Suppose that $M$ is a totally orientable graph manifold and $M = V
\cup_S W$ is a Heegaard splitting.  Suppose further that there is an
exceptional fiber $f$ in $M_v$ and an annulus $A$ such that:

{\rm(1)}\qua One component of $\partial A$ wraps at least twice around $f$; and

{\rm(2)}\qua $A$ is embedded away from $\partial A \cap f$; and

{\rm(3)}\qua $\partial A \backslash f$ lies in $S$.  

Then $f$ is a core of either $V$ or $W$.
\end{lem}

\proof In fact, in the central argument in \cite{Sch1}, the existence
of a boundary component is used exclusively to produce such an
annulus.
\endproof

This more general lemma will be used in the next subsection.  A
consequence of Lemma \ref{lem:bcore} is that we can make use of the
following lemma.

\begin{lem} \label{lem:ex}
Suppose that $f$ is an exceptional fiber of $M_v$ and that $f$ is also
a core of $V$.  Then $M \backslash \eta(f) = (V \backslash \eta(f))
\cup_S W$ is a Heegaard splitting.  Furthermore, $S \cap (M_v
\backslash \eta(f))$ is vertical or pseudovertical, respectively, if
and only if $S \cap M_v$ is vertical or pseudovertical, respectively.
The same holds if $f$ is a core of $W$.
\end{lem}

\proof Since $f$ is a core of $V$, $V \backslash \eta(f)$ is still a
compression body.  Thus $M \backslash \eta(f) = (V \backslash \eta(f))
\cup_S W$ is a Heegaard splitting.  Conversely, if $M \backslash
\eta(f) = (V \backslash \eta(f)) \cup_S W$ is a Heegaard splitting,
then $(V \backslash \eta(f)) \cup N(f)$ is a compression body.  Hence
$M = V \cup_S W$ is a Heegaard splitting.

\begin{figure}[ht!]
\centerline{\epsfxsize=1.8in \epsfbox{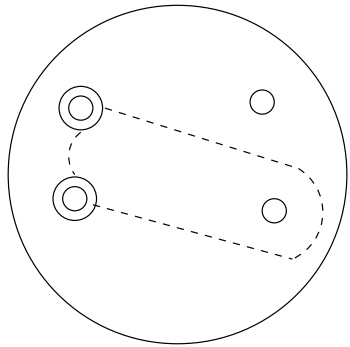}}
\caption{Schematic for vertical Heegaard splittings of graph manifold
with boundary}
\label{gm11b.ps}
\end{figure}

Now by the definition of vertical and pseudovertical, respectively, $S
\cap (M_v \backslash \eta(f))$ is vertical or pseudovertical,
respectively, if and only if $S \cap M_v$ is vertical or
pseudovertical, respectively.  Compare Figures \ref{gm11a.ps} and
\ref{gm11b.ps}.  \endproof

The following Proposition generalizes Theorem \ref{thm:sch1}.  But
because Theorem \ref{thm:sch1} is already known, we restrict our
attention to the case in which $M$ has a non empty characteristic
submanifold.

\begin{prop} \label{prop:eb}
Suppose that $M$ is a connected generalized graph manifold with non
empty characteristic submanifold.  Suppose that $M = V \cup_S W$ is a
strongly irreducible Heegaard splitting.  Suppose that the active
component of $M = V \cup_S W$ is a vertex manifold $M_v$ that is a
Seifert fibered space and that has exterior boundary.  Then we may
isotope the decomposing tori, thereby redefining $M_v$ and the edge
manifolds for which $e$ is incident to $v$ slightly, so that after
this isotopy, $S \cap M_v$ is a vertical surface and so that an edge
manifold becomes the active component.
\end{prop}

\proof The proof is by induction on the number of exceptional fibers
in $M_v$.  Suppose first that $M_v$ contains no exceptional fibers.
Denote the exterior boundary of $M_v$ by $\partial_E M_v$.  Let $\cal
A$ be a collection of disjoint essential vertical annuli in $M_v$ that
cut $M_v$ into a regular neighborhood of $\partial M_v \backslash
\partial_E M_v$.  After an isotopy, $S \cap \cal A$ consists of closed
curves essential in both $S$ and $\cal A$ and in the minimal possible
number of such curves.  In particular, after a small isotopy, this
intersection consists of regular fibers of $M_v$.

\begin{figure}[ht!]
\centerline{\epsfxsize=1.8in \epsfbox{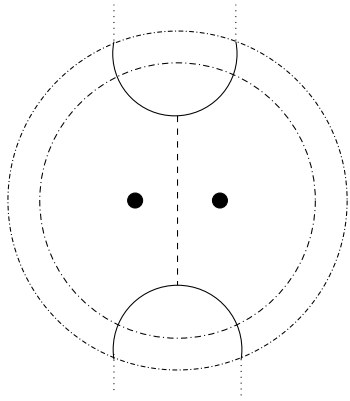}}
\caption{The active component is a vertex manifold}
\label{gm10b.ps}
\end{figure}

Now isotope $S$ so that $S \cap N(\cal A)$ consists of vertical
incompressible annuli.  Set $\tilde M_v = N(\partial_E M_v) \cup
N(\cal A)$.  Isotope $S$ so that $S \cap \tilde M_v$ consists of $S
\cap N(\cal A)$ together with annuli in $N(\partial_E M_v) \backslash
N(\cal A)$ that join two components of $S \cap N(\cal A)$.  Isotope
any annuli in $S \cap (M_v \backslash \tilde M_v)$ that are parallel
into $\partial \tilde M_v$ into $\tilde M_v$.  Set $\partial_E \tilde
M_v = \partial_E M_v$.

\begin{figure}[ht!]
\centerline{\epsfxsize=1.8in \epsfbox{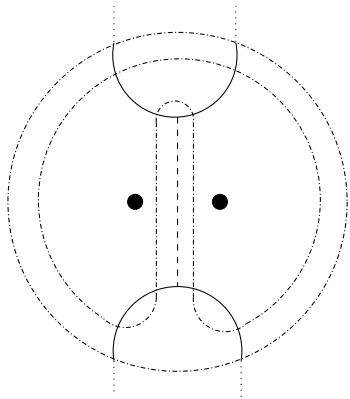}}
\caption{The active component is an edge manifold}
\label{gm10c.ps}
\end{figure}

Let $\tilde T$ be a component of $\partial \tilde M_v \backslash
\partial_E \tilde M_v$.  Then $\tilde T$ is parallel to a decomposing
torus or annulus $T$.  We replace $T$ by $\tilde T$.  We may do so via
an isotopy.  We do this for all components of $\partial \tilde M_v
\backslash \partial_E \tilde M_v$.  After this process, the
conclusions of the proposition hold.  See Figures \ref{gm10b.ps} and 
\ref{gm10c.ps}.

To prove the inductive step, suppose that $f$ is an exceptional fiber
of $M_v$.  Then by Lemma \ref{lem:bcore}, $f$ is a core of either $V$
or $W$, say of $V$.  The inductive hypothesis in conjunction with
Lemma \ref{lem:ex} then proves the theorem.  \endproof

\subsection{Seifert fibered vertex manifold without exterior boundary}

Next we consider the case in which the active component of $M = V
\cup_S W$ is a vertex manifold $M_v$ that is a Seifert fibered space
and that has no exterior boundary.  Part of the strategy is somewhat
reminiscent of the strategy used in the preceding case.  However, the
situation here is more complicated and a more refined strategy must be
used.  The refined strategy involves a generalization of Lemma
\ref{lem:essential} to a ``spine'' of the vertex manifold.

\begin{defn}
Let $Q_v$ be the base orbifold of $M_v$.  A spine for $Q_v$ is a
1--complex $\Gamma^1$ with exactly one vertex $v$ that cuts $Q_v$ into
a regular neighborhood of $\partial Q_v \cup (exceptional \; points)$.
A 2--complex of the form $\Gamma^2 = p^{-1}(\Gamma^1)$ is called a
spine of $M_v$.  We denote the regular fiber $p^{-1}(v)$ by $\gamma$.
\end{defn}

The following lemma generalizes Lemma \ref{lem:essential}.  Though we
will only be interested in this lemma in the case that $N$ is a vertex
manifold of a graph manifold, we state it in very general terms.  It
applies to a larger class of 3--manifolds than just graph manifolds.

\begin{lem} \label{lem:gessential}
Let $M = V \cup_S W$ be a strongly irreducible Heegaard splitting.
Let $N$ be a totally orientable Seifert fibered submanifold of $M$
that doesn't meet $\partial M$.  Let $\Gamma^2$ be a spine of $N$.
Then $S$ may be isotoped so that the following hold:

{\rm(1)}\qua $S \cap \Gamma^2$ consists of simple
closed curves and simple closed curves wedged together at points in
$\gamma$.

{\rm(2)}\qua No closed curve in $S \cap \Gamma^2$ bounds a disk in $\Gamma^2
   \backslash S$.
\end{lem}

\proof 
The first part of the assertion follows by general position.
To prove the second assertion, let $X$ be a spine of $V$ and $Y$ a
spine of $W$.  Then $M \backslash (\partial_-V \cup X \cup \partial_-W
\cup Y)$ is homeomorphic to $S \times (0, 1)$.  $X$ can't be disjoint
from $\Gamma^2$.  Thus for $t$ near $0$, $(S \times t) \cap \Gamma^2$
contains simple closed curves that bound essential disks in $V$.  $Y$
can't be disjoint from $\Gamma^2$ either.  Thus for $t$ near $1$, $(S
\times t) \cap \Gamma^2$ contains simple closed curves that bound
essential disks in $W$.  

As $t$ increases, $(S \times t) \cap \Gamma^2$ changes continuously.
Since $M = V \cup_S W$ is strongly irreducible, there can be no $t$
such that $(S \times t) \cap \Gamma^2$ contains simple closed curves
that bound essential disks in $V$ and simple closed curves that bound
essential disks in $W$.  Thus, there is a $t_0$, such that $(S \times
t_0) \cap \Gamma^2$ contains no simple closed curves that bound
essential disks in $V$ or $W$.  Any remaining disk components in $(S
\times t_0) \cap \Gamma^2$ must be inessential in $V$ or $W$ and can
hence be removed via isotopy.  The lemma follows.  \endproof

Before launching into the two main portions of the argument, we prove
an auxiliary lemma.  This lemma is a weak version of a counterpart to
Theorem 3.3 in \cite{Sc}.

\begin{figure}[ht!]
\centerline{\epsfxsize=2.0in \epsfbox{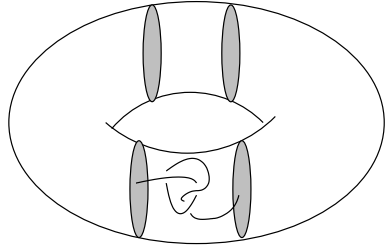}}
\caption{Schematic for a Heegaard splitting intersecting a solid
torus}
\label{gm15.ps}
\end{figure}

\begin{lem} \label{lem:solidtorus}
Suppose $M = V \cup_S W$ is a strongly irreducible Heegaard splitting
of a 3--manifold $M$.  Suppose $U \subset M$ is a solid torus such that
$S$ intersects $\partial U$ in meridians.  Further suppose that
$\partial (M \backslash interior(U))$ is incompressible in $M
\backslash interior(U)$.  Then $S \cap U$ consists of meridian disks
of $U$ and components that are obtained by ambient 1--surgery on pairs
of meridian disks along a single arc that joins the two meridian disks.
\end{lem}

\proof
Let $s$ be a component of $S \cap \partial U$. A collar of $s$ in $S$
is an annulus $A$.  Lemma 2.6 in \cite{ST1} states: ``Suppose $S$
gives a Heegaard splitting of a 3--manifold $M$ into compression bodies
$V$ and $W$.  Suppose that $F \subset S$ is a compact subsurface so
that every component of $\partial F$ is essential in $S$.  Suppose
each component of $\partial F$ bounds a disk in $M$ disjoint from
$interior(F)$.  Either $\partial F$ bounds a collection of disks in a
single compression body or $M = V \cup_S W$ is weakly reducible.''

Here $M = V \cup_S W$ is strongly irreducible.  It follows that either
$s$ bounds a disk in $S$ or $s$ bounds a disk in a single compression
body, say $V$.  Here $\partial (M \backslash interior(U))$ is
incompressible in $M \backslash interior(U)$.  It follows that in case
of the former, $s$ bounds a meridian disk of $U$ in $S \cap U$ and
that in case of the latter, a meridian disk of $U$ bounded by $s$ lies
entirely in $V$.

Thus either $s$ bounds a meridian disk in $S \cap U$ or bounds a
meridian disk in the surface obtained by compressing $S \cap U$ along
$D$.  It follows that $S \cap U$ may be reconstructed from meridian
disks by ambient 1--surgery along a collection of arcs such that each
meridian disk meets at most one endpoint of one of the arcs.  In
particular, between two meridian disks there is at most one arc.
\endproof

\begin{lem} 
Suppose that the active component of $M = V \cup_S W$ is a vertex
manifold $M_v$ that is a Seifert fibered space with boundary but 
no exterior boundary.  Then one of the following holds:

{\rm(1)}\qua We may isotope the decomposing tori, thereby redefining $M_v$ and
the edge manifolds for which $e$ is incident to $v$ slightly, so that
after this isotopy, $S \cap M_v$ is either a horizontal or vertical
incompressible surface and so that an edge manifold becomes the active
component.

Or:

{\rm(2)}\qua $S$ may be isotoped within $M_v$ so that a fiber $f$ of $M_v$ lies
   in $S$.

\end{lem}

\proof Let $\Gamma^2$ be a spine of $M_v$ and isotope $S$ within $M_v$
so that the conclusions of Lemma \ref{lem:gessential} hold.  Three
cases need to be considered:

{\bf Case 1}\qua A simple closed curve in $S \cap \Gamma^2$ can be isotoped to
be vertical.

Then $S$ may be isotoped in $M_v$ so that it contains a
regular fiber of $M_v$.

{\bf Case 2}\qua No simple closed curve in $S \cap \Gamma^2$ can be isotoped to
be vertical and in each solid torus component $U$ of $M_v \backslash
\eta(\Gamma^2)$, $S \cap \partial U$ consists of meridians.

In this case, after a small isotopy, $S \cap
N(\Gamma^2)$ is a horizontal incompressible surface.  Furthermore, let
$U$ be a component of $M_v \backslash \eta(\Gamma^2)$.  By Lemma
\ref{lem:solidtorus} $S \cap U$ consists of meridian disks possibly
together with other components that are obtained by ambient 1--surgery
on pairs of meridian disks along a single arc that joins the two
meridian disks.  Isotope each such arc out of $U$, through
$N(\Gamma^2)$, avoiding $\gamma$, to lie in a component of $N(\partial
M_v)$.  Then $S \cap U$ consists of meridians for each solid torus
component of $M_v \backslash \Gamma^2$ and all ambient 1--surgeries
occur in $N(\partial M_v)$.

Let $\tilde M_v$ be the union of $N(\Gamma^2)$ with the solid tori
containing the exceptional fibers of $M_v$.  Then $\tilde M_v$ is a
shrunk version of $M_v$.  Note that $S \cap \tilde M_v$ is a
horizontal incompressible surface.  Let $\tilde T$ be a component of
$\partial \tilde M_v$.  Then $\tilde T$ is parallel to a decomposing
torus $T$.  We replace $T$ by $\tilde T$.  We may do so via an
isotopy.  We do this for all components of $\partial \tilde M_v$.
After this process, the conclusions of the lemma hold.

{\bf Case 3}\qua No simple closed curve in $S \cap \Gamma^2$ can be isotoped to
be vertical and there is a solid torus component $U$ of $M_v
\backslash \eta(\Gamma^2)$ such that $S \cap \partial U$ does not
consist of meridians.

The argument in this case is given, for instance, in
Proposition 1.1 of \cite{BO1}.  For completeness we povide a sketch of
an argument more in line with the ideas used here.  In this case there
is a possibly singular annulus $A$ between a component of $S \cap
\partial U$ and the exceptional fiber $f$ in $U$.  If $A$ is not
singular, then $A$ describes an isotopy of $S$ after which $f$ lies in
$S$.  If $A$ is singular, then $A$ satisfies the hypotheses of Lemma
\ref{lem:core}.  Thus $f$ is a core of either $V$ or $W$.  Ie, $S$
is the splitting surface for a Heegaard splitting of $M \backslash
\eta(f)$.  But then Proposition \ref{prop:eb} applies to $M_v
\backslash \eta(f)$ and we may isotope the decomposing tori, thereby
redefining $M_v$ and the edge manifolds for which
$e$ is incident to $v$ slightly, so that after this isotopy, $S \cap
M_v$ is a vertical incompressible surface and so that an edge manifold
becomes the active component.
\endproof

\begin{lem} \label{lem:neb}
Suppose that the active component of $M = V \cup_S W$ is a vertex
manifold $M_v$ that is a Seifert fibered space with boundary but no
exterior boundary.  Suppose further that a fiber $f$ of $M_v$ lies in
$S$.  Then one of the following holds:

{\rm(1)}\qua We may isotope the decomposing tori, thereby redefining the active
component $M_v$ and the edge manifolds for which $e$ is incident to
$v$ slightly, so that after this isotopy, $S \cap M_v$ is a vertical
incompressible surface and so that an edge manifold becomes the active
component.

Or:

{\rm(2)}\qua $S \cap M_v$ is pseudohorizontal.  
\end{lem}

\proof Consider a small regular neighborhood $N(f)$ of $f$ such that
$S \cap N(f)$ is a collar $A$ of $f$.  Compress $S$ as much as
possible in $M_v \backslash \eta(f)$ to obtain an incompressible
surface $S^* \subset M \backslash \eta(f)$.  By Haken's Theorem, each
compressing disk can be chosen to lie entirely on one side of $S$.
Since $M = V \cup_S W$ is strongly irreducible, all compressions must
have been performed to one side of $S$.  It follows that $S^*$ lies
either in $V$ or in $W$.  There are three options for $S^* \cap M_v$:

{\bf Case 1}\qua $S^* \cap (M_v \backslash \eta(f))$ contains an annulus that
is parallel into $\partial N(f)$.

Note that $S \cap N(f)$ is also parallel into $\partial N(f)$.  Thus
$S \cap M_v = (S^* \cap (M_v \backslash \eta(f))) \cup (S \cap N(f))$
contains a torus bounding a solid torus in either $V$ or $W$.
Furthermore, $f$ lies on the boundary of this solid torus and meets a
meridian disk once.  After a small isotopy, $f$ is a core of the solid
torus.  Thus $f$ is a core of either $V$ or $W$.  But then Proposition
\ref{prop:eb} applies to $M_v \backslash \eta(f)$ and we may isotope
the decomposing tori, thereby redefining $M_v$ and the edge manifolds
for which $e$ is incident to $v$ slightly, so that after this isotopy,
$S \cap M_v$ is vertical and so that an edge manifold becomes the
active component.

{\bf Case 2}\qua $S^* \cap (M_v \backslash \eta(f))$ is vertical.

If $S^* \cap (M_v \backslash \eta(f))$ is not boundary parallel, then
it is essential and can't be contained in $V$ or $W$.  This is a
contradiction, hence this case does not occur.

{\bf Case 3}\qua $S^* \cap (M_v \backslash \eta(f))$ is horizontal.

Then $(S^* \cap (M_v \backslash \eta(f))) \cup (S \cap N(f))$ is
pseudohorizontal, but we must show that in fact $S \cap M_v$ $= (S^*
\cap (M_v \backslash \eta(f))) \cup (S \cap N(f))$.

Since $\partial A$ has two components
and since $S^* \cap (M_v \backslash \eta(f))$ is separating there are two,
necessarily parallel, components of $S^* \cap (M_v \backslash \eta(f))$.
Moreover, since $A$ is parallel into $\partial U$ in both directions,
each component of $M_v \backslash S^*$ is homeomorphic to $(punctured
\; \surface) \times (0, 1)$.

Suppose now that $S^*$ lies in, say, $V$.  Denote the component of $M
\backslash S^*$ that meets $W$ by $\hat W$.  Then $S$ defines a
Heegaard splitting of $\hat W$.  Let ${\cal D}$ be a set of disks in
the interior of $M_v$ that cut $\hat W \cap M_v$ into a collar of the
annuli $\hat W \cap \partial M_v$.  By Haken's Theorem, each such disk
can be isotoped to intersect $S$ in a single circle.  We may assume
that after this isotopy, ${\cal D}$ still lies in the interior of
$M_v$.  

Here $S$ may be reconstructed by performing ambient 1--surgery on $S^*$
along arcs in $\hat W$.  But this collection of arcs is disjoint from
${\cal D}$.  Thus all such arcs may be isotoped into edge manifolds
$M_e$ such that $e$ is incident to $v$.  Note that an edge manifold
that contains such an arc becomes the active component.  Also note
that after this isotopy, the portion of $S$ remaining in $M_v$ is
pseudohorizontal.

If a surface is pseudohorizontal, then it is boundary compressible.
In particular, it lives in the active component.  Hence the existence
of such arcs contradicts Lemma \ref{lem:C}.  Thus $S \cap M_v$ $= (S^*
\cap (M_v \backslash \eta(f))) \cup (S \cap N(f))$.  \endproof

\subsection{Vertex manifold not Seifert fibered}

Next we consider the case in which the active component of $M = V
\cup_S W$ is a vertex manifold $M_v$ that is homeomorphic to $(compact
\; orientable \; \surface) \times I$.  The strategy here is an
adaptation of the argument in the preceeding case.  Though the setup
here is much simpler.

\begin{defn}
Let $Q$ be a compact surface with non empty boundary.  A spine of $Q$
is a 1--complex $\Gamma^1$ with exactly one vertex $v$ that cuts $Q$
into a regular neighborhood of $\partial Q$.  A 2--complex of the form
$\Gamma^2 = \Gamma^1 \times I$ is called a spine of $Q \times I$.  We
denote the 1--manifold $v \times I$ by $\gamma$.
\end{defn}

The following lemma is another generalization of Lemma
\ref{lem:essential}.  We will only be interested in this lemma in the
case that $N$ is a vertex manifold of a graph manifold, but again we
state it in very general terms.  It too applies to a larger class of
3--manifolds than just generalized graph manifolds.

\begin{lem} \label{lem:ggessential}
Let $M = V \cup_S W$ be a strongly irreducible Heegaard splitting.
Suppose $N$ is a submanifold of $M$ homeomorphic to $(compact\; \surface)
\times I$ such that $(compact \; \surface) \times \partial I \subset
\partial M$.  Let $\Gamma^2$ be a spine of $N$.  
Then $S$ may be isotoped so that the following hold:

{\rm(1)}\qua $S \cap \Gamma^2$ consists of simple
closed curves and simple closed curves wedged together at points in
$\gamma$.

{\rm(2)}\qua No closed curve in $S \cap \Gamma^2$ bounds a disk in $\Gamma^2
   \backslash S$.
\end{lem}

The proof of this lemma is identical to the proof of Lemma
\ref{lem:gessential}.

\begin{lem} \label{lem:prod}
Suppose that the active component of $M = V \cup_S W$ is a vertex
manifold $M_v$ that is homeomorphic to $(compact \; orientable \;
\surface) \times I$.  Then we may isotope the decomposing annuli,
thereby redefining the active component $M_v$ and the edge manifolds
for which $e$ is incident to $v$ slightly, so that after this isotopy,
$S \cap M_v$ is a horizontal incompressible surface and so that an
edge manifold becomes the active component.
\end{lem}

\proof 
Let $\Gamma^2$ be a spine of $M_v$ and isotope $S$ within $M_v$
so that the conclusions of Lemma \ref{lem:ggessential} hold.  Then $S
\cap \Gamma^2$ consists of horizontal curves.  Here $S \cap
N(\Gamma^2)$ is a bicollar of $S \cap \Gamma^2$ and hence a horizontal
incompressible surface.  

Let $\tilde M_v$ be $N(\Gamma^2)$.  Let $\tilde A$ be a component of
$\partial \tilde M_v \backslash \partial M_v$.  Then $\tilde A$ is
parallel to a decomposing annulus $A$.  We replace $A$ by $\tilde A$.
We may do so via an isotopy.  We do this for all components of
$\partial \tilde M_v \backslash \partial M_v$.  After this process,
the conclusions of the lemma hold.
\endproof

\subsection{Edge manifold homeomorphic to $(annulus) \times I$}

Now we consider the case in which the active component of $M = V
\cup_S W$ is an edge manifold $M_e$ homeomorphic to $(annulus) \times
I$.  This case turns out to be a direct application of the following
theorem of Marty Scharlemann:

\begin{thm} \label{thm:Sc}
Suppose $M = V \cup_S W$ is a strongly irreducible Heegaard splitting
and $U \subset M$ is a solid torus such that $S$ intersects $\partial
U$ in parallel essential non meridional curves.  Then $S$ intersects
$U$ in a collection of boundary parallel annuli and possibly one other
component, obtained from one or two annuli by ambient 1--surgery along
an arc parallel to a subarc of $\partial U$.  If the latter sort of
component is in $U$, then $S \backslash U$ is incompressible in $M
\backslash U$.
\end{thm}

\proof
This is Theorem 3.3 in \cite{Sc}.
\endproof

The following proposition is stated in general terms.  Our interest in
this theorem will be the case in which $A \times I$ is an edge
manifold.  Since $A \times I$ is a solid torus, this proposition
follows directly from Scharlemann's Theorem (Theorem \ref{thm:Sc}).

\begin{prop} \label{prop:acs}
Suppose $A \times I$ is an imbedding of $(annulus) \times I$ in the
interior of $M$ with $A \times \{point\}$ essential in $M$.  Further
suppose that $M = V \cup_S W$ is a strongly irreducible Heegaard
splitting.  If both components of $S \cap (A \times \partial I)$
consist of curves essential in both $S$ and $(A \times \partial I)$,
then $S \cap (A \times I)$ is isotopic to a collection of
incompressible annuli and possibly one other component, obtained from
two annuli by ambient 1--surgery along an arc parallel to a
subarc of $A \times \{point\}$.
\end{prop}

The following definition clarifies the statement of Proposition
\ref{prop:acs}.

\begin{defn}
Let $A$ be an annulus and let $N = A \times I$.  A spanning annulus
for $N$ is an annulus of the form $c \times I$, for $c$ an essential
curve in $A$.  Similarly, let $T$ be a torus and $N = T \times I$.  A
spanning annulus for $N$ is an annulus of the form $c \times I$, for
$c$ an essential curve in $T$.
\end{defn}

\begin{figure}[ht!]
\centerline{\epsfxsize=5.0in \epsfbox{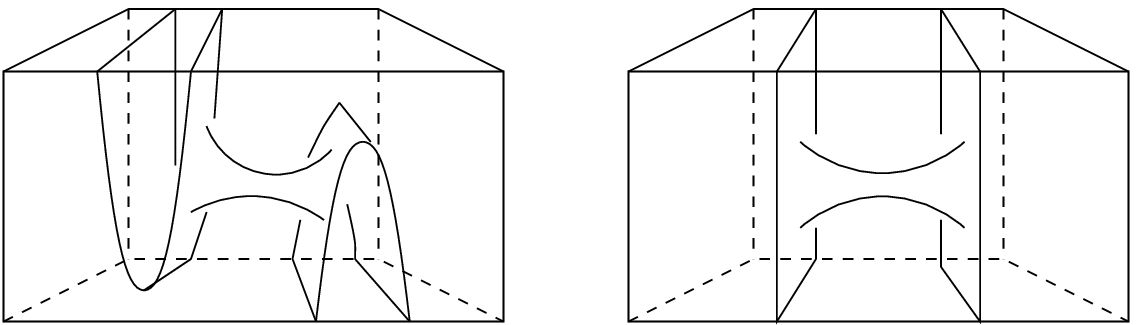}}
\caption{Two isotopic possibilities}
\label{gmx2.ps}
\end{figure}

Each of the incompressible annuli mentioned in Proposition
\ref{prop:acs} is either a spanning annulus or is parallel into $A
\times \partial I$.  In Figure \ref{gmx2.ps} we see two isotopic
possibilities.  Note that the tube we see on the left hand side is
actually ``dual'' to the tube we see on the right hand side.  

The tube on the left hand side could (after straightening out the
picture) be seen as ambient 1--surgery along an arc parallel to a
subarc of $\{point\} \times I$ (a vertical arc).  In this case the the
tube on the right hand side would be seen as ambient 1--surgery along
an arc parallel to a subarc of $A \times \{point\}$ (a horizontal arc).

\subsection{Edge manifold homeomorphic to $(torus) \times I$}

Next we consider the case in which the active component of $M = V
\cup_S W$ is an edge manifold $M_e$ homeomorphic to $(torus) \times
I$.  The theorem proven here applies in very general contexts.  The
techniques used are those developed by Rubinstein and Scharlemann in
\cite{RS}.  

The use of these techniques is inspired by Cooper and Scharlemann's
application of the central argument in \cite{RS} to the setting of
Heegaard splittings of solvmanifolds in \cite[Theorem 4.2]{CS}.
The arguments in this section are identical to the argument in
\cite[Theorem 4.2]{CS}.  We weaken one of the hypotheses slightly.  In
addition, in the setting here, there are no constraints on $S \cap (T
\times \partial I)$.  In \cite{CS} a constraint arises due to the fact
that in a solvmanifold the two components of $T \times \partial I$ are
identified.  This means that some scenarios that arise in the argument
for \cite[Theorem 4.2]{CS} can be ruled out there.  But they can't be
ruled out here.  Thus our conclusions are slightly different.

The argument is rather lengthy.  The reader is referred to \cite{CS}
for a sketch and to \cite{RS} for details.  

We recall some fundamentals concerning the Rubinstein--Scharlemann
graphic.  Let $\Sigma_V$ be a spine of $V$ and $\Sigma_W$ a spine of
$W$.  Then $M \backslash (\partial_-V \cup \Sigma_V \cup \partial_-W \cup
\Sigma_W)$ is homeomorphic to $S \times I$.  This product foliation is
called a \emph{sweepout}.  If $S$ intersects a product
submanifold $N = Q \times I$ of $M$ then to each point $(s, t)$ in $I
\times I$ we may associate $S_s = S \times \{s\}$ and $Q_t = Q \times
\{t\}$.  We are interested in $S_s \cap Q_t$.

\begin{figure}[ht!]
\centerline{\epsfxsize=1.5in \epsfbox{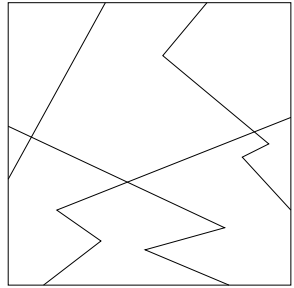}}
\caption{The Rubinstein--Scharlemann graphic}
\label{gma.eps}
\end{figure}

After a small isotopy, we may assume that the spines and the sweepout
are in general position with respect to the foliation of $Q \times I$.
There is then a $1$--dimensional complex $\Gamma$ in $I \times I$
called the \emph{Rubinstein--Scharlemann graphic}. See Figure
\ref{gma.eps}.  At a point $(s, t)$ away from $\Gamma$, $S_s$ and
$Q_t$ are in general position.  On an edge of $\Gamma$ the surfaces
$S_s$ and $Q_t$ have a single point of tangency.  At a vertex of
$\Gamma$ there are either two points of tangency or a ``birth--death''
singularity.

For $(s,t)$ in a region of $(I \times I) \backslash \Gamma$, $S_s \cap
Q_t$ is topologically rigid.  If there is a curve $c$ in $S_s \cap
Q_t$ that is essential in $S_s$ but bounds a disk in $Q_t$ that lies
in $V$ near $c$, then we label the region $V$.  Similarly, if there is
a curve $c'$ in $S_s \cap Q_t$ that is essential in $S_s$ but bounds a
disk in $Q_t$ that lies in $W$ near $c$, then we label the region $W$.

We summarize some of the insights from \cite{RS}:

\begin{rem} \label{rem:labeldisk}
If a region $R$ is labelled $V$, then for $(s, t) \in R$ the
curve $c$ in $S_s \cap Q_t$ that is essential in $S_s$ and bounds a
disk in $Q_t$ that lies in $V$ near $c$ also bounds a disk that
lies entirely in $V$.  The equivalent statement holds for the label
$W$.  (This is \cite[Lemma 4.3]{RS}.)
\end{rem}

This fact has the following immediate consequences.

\begin{rem} \label{rem:notvw}
If a region is labelled both $V$ and $W$, then $M = V \cup_S W$ is
weakly reducible.  (This is \cite[Corollary 4.4]{RS}.)
\end{rem}

\begin{rem} \label{rem:mon}
If for some $s_0 \in (0, 1)$ and some $t_0, t_1 \in (0, 1)$ there is a
region labelled $V$ containing $(s_0, t_0)$ and a region labelled $W$
containing $(s_0, t_1)$, then $M = V \cup_S W$ is weakly reducible.
\end{rem}

Remark \ref{rem:mon} is not explicitly stated in \cite{RS} but follows
immediately from the definitions and Remark \ref{rem:labeldisk}.
Indeed, the labelling in the region containing $(s_0, t_0)$ gives an
essential disk to one side of $S_{s_0}$ and the labelling in the
region containing $(s_0, t_1)$ gives an essential disk to the other
side of $S_{s_0}$.  The boundaries of these disks are contained in
$Q_{t_0}$ and $Q_{t_1}$ respectively.  As $Q_{t_0}$ and $Q_{t_1}$ are
disjoint, the boundaries of these disks are disjoint.  Hence $M = V
\cup_S W$ is weakly reducible.

More subtle facts are the following:

\begin{rem} \label{rem:adj}
The labels $V$ and $W$ cannot both appear in
regions adjacent along an edge.  (This is \cite[Corollary 5.5]{RS}.)
\end{rem}

Notice that for small $s$, $S_s \cap (Q \times I)$ is the boundary of
a small regular neighborhood of $\Sigma_V \cap (Q \times I)$.  If the
leaves of the foliation of $Q \times I$ are essential, then $\Sigma_V$
can't miss any such leaf.  Similarly for $\Sigma_W$.  Thus we have the
following:

\begin{rem} \label{rem:vw}
If the leaves of the foliation of $Q \times I$ are essential in $M$,
then every region of $(I \times I) \backslash \Gamma$ abutting the
left edge of $I \times I$ is labelled $V$.  Every region of $(I \times
I) \backslash \Gamma$ abutting the right edge of $I \times I$ is
labelled $W$.
\end{rem}

As a warm up in the application of these principles, we prove the
following lemma:

\begin{lem} \label{lem:vertex}
Suppose that $N = Q \times I$ is a product submanifold of $M$.
Suppose further that $M = V \cup_S W$ is a strongly irreducible
Heegaard splitting.  Let $\Gamma$ be the Rubinstein--Scharlemann
graphic and suppose there is a vertex $(s_0, t_0)$ of $\Gamma$
with the following properties:

{\rm(1)}\qua Four regions meet at $(s_0, t_0)$.  

{\rm(2)}\qua One region is labelled $V$ and one region is labelled $W$.  

Then the following hold:

{\rm(A)}\qua The other two regions abutting $(s_0, t_0)$
are unlabelled.

{\rm(B)}\qua The regions labelled $V$ and $W$ lie opposite each other.

{\rm(C)}\qua The graph $G = S_{s_0} \cap Q_{t_0}$ in $Q_{t_0}$ contains a
connected subgraph $\tilde G$ with two vertices $v_1, v_2$, each of
valence 4.

{\rm(D)}\qua If an edge $e$ of $G$ has both ends on the same vertex $v$ then
$e \cup v$ is an essential circle in $Q_{t_0}$.

{\rm(E)}\qua $\tilde G$ is not contractible in $Q_{t_0}$.
\end{lem}

\proof Observations A and B follow directly from Remark \ref{rem:adj}.
Recall that each edge signifies a tangency between $S_{s_0}$ and
$Q_{t_0}$.  If the tangency corresponds to a maximum or minimum, then
the labelling of adjacent regions does not change as the edge is
traversed.  So here the tangencies corresponding to edges must arise
from saddle singularities.

\begin{figure}[ht!]
\centerline{\epsfxsize=1.5in \epsfbox{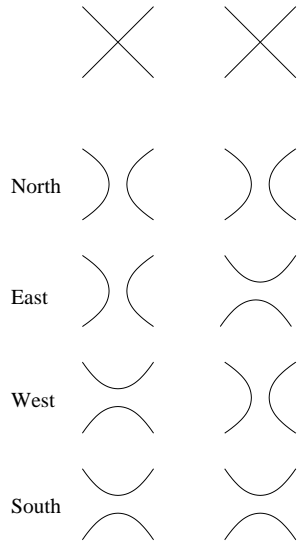}}
\caption{How $S_s$ meets $T_t$ in the four regions North, East,
West and South meeting $(s_0, t_0)$}
\label{gm1.eps}
\end{figure}

It follows that $S_{s_0} \cap Q_{t_0}$ consists of a graph $G$ in
$Q_{t_0}$ that has two valence four vertices.  In the four regions
abutting $(s_0, t_0)$ the valence four vertices break apart in four
possible combinations.  See Figure \ref{gm1.eps}.

Suppose the two valence four vertices of $G$ lie on distinct
components $\tilde G_1, \tilde G_2$ of $G$.  As $\tilde G_1, \tilde
G_2$ break apart as in Figure \ref{gm1.eps}, the curves arising from
$\tilde G_1$ remain disjoint from the curves arising from $\tilde
G_2$.  Furthermore, the curves arising from $\tilde G_1$ are identical
in the regions North and East and in the regions
West and South.  The curves arising from $\tilde G_2$ are identical in
the regions North and West and in the regions East
and South.  It follows that one of the components, say $\tilde G_1$,
must give rise to the labelling $V$ and the other component, say
$\tilde G_2$, must give rise to the labelling $W$.  But then two
adjacent regions, say the regions North and East, are
labelled $V$.  A contradiction.  Thus C holds.

Denote the vertices by $v_1$ and $v_2$.  Suppose that $e$ has both
ends on $v_1$ and that $v_1 \cup e$ is inessential.  Then there will
be two adjacent regions, say the regions North and East, in which this
monogon gives rise to a simple closed curve.  Since one of the
regions, say North, is labelled, say $V$, this simple closed curve
must give rise to a labelling $V$.  But then it gives rise to this
labelling in both the region North and in the region East.  A
contradiction.  Hence D holds.

\begin{figure}[ht!]
\centerline{\epsfxsize=1.5in \epsfbox{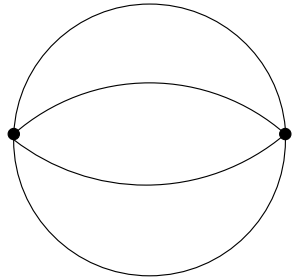}}
\caption{$\tilde G$ is contractible}
\label{gm1a.eps}
\end{figure}

\begin{figure}[ht!]
\centerline{\epsfxsize=5.0in \epsfbox{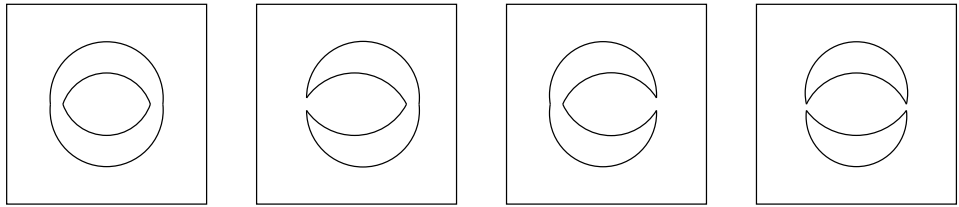}}
\caption{The regions North, East, West and South}
\label{gm1b.eps}
\end{figure}

Finally, suppose $\tilde G$ is contractible in $Q_{t_0}$.  By D, each
edge has endpoints on distinct vertices.  Thus $\tilde G$ is as in
Figure \ref{gm1a.eps}.  Hence in the regions North, East, West and
South of $(I \times I) \backslash \Gamma$ we see components of
intersection of $S_s \cap Q_t$ in $Q_t$ as in Figure \ref{gm1b.eps}.

The curves of intersection pictured are the ones giving rise to the
labellings.  Note that if the circles pictured in East and West give
rise to labellings, then they give rise to the same labelling.  Hence
these regions must both be unlabelled.  It follows that North and
South are the labelled regions.  If the larger circle pictured in
North gives rise to a labelling, then this labelling coincides with
the labelling arising from the circle(s) pictured in South.  But this is
impossible.  Thus the smaller circle pictured in North gives rise to
the labelling.  We may assume that it gives rise to the label $V$.

Let $(s, t) \in$ North.  Denote the annulus cut out by the two circles
in North by $\tilde A$.  Denote the disk cut out by the inner circle
by $\tilde D$.  By Remark \ref{rem:labeldisk} we may assume that
$\tilde D \subset V$.  Denote the outer component of $\partial \tilde
A$ by $c$ and the inner component by $c'$.  Since $c'$ gives rise to a
labelling, $c'$ is essential in $S_s$.  Since $c$ does not give rise
to a labelling, $c$ is inessential in $S_s$.  Let $D \subset S_s$ be
the disk bounded by $c$.  We may assume that this disk is disjoint
from $c \cup c'$.  After a small isotopy, $\tilde A \cup D$ is an
essential disk in $W$.  But here $\partial (\tilde A \cup D) =
\partial \tilde D$.  So $M = V \cup_S W$ is reducible.  But this is
impossible.  Thus $\tilde G$ is not contractible in $Q_{t_0}$.  Hence
E holds.
\endproof

The following proposition gives the required information concerning
the intersection of a strongly irreducible Heegaard splitting with an
edge manifold homeomorphic to $(torus) \times I$ when this edge
manifold is the active component.  

\begin{prop} \label{prop:tcs}
Suppose $T \times I$ is an imbedding of $(torus) \times I$ in the
interior of $M$ with $T \times \{t\}$ essential in $M$.  Further
suppose that $M = V \cup_S W$ is a strongly irreducible Heegaard
splitting.  If $S \cap (T \times \partial I)$ consist of curves
essential in both $S$ and $(T \times \partial I)$, then $S \cap (T
\times I)$ is characterized by one of the following:

{\rm(1)}\qua $S \cap (T \times I)$ is isotopic to a collection of incompressible
annuli and possibly one other component, obtained from two such annuli
by ambient 1--surgery along an arc parallel to a subarc of $T \times
\{point\}$.

{\rm(2)}\qua There is a pair of simple closed curves $c, c' \subset T$ such that
   $c \cap c'$ consists of a single point $p \in T$ and $V \cap (T
   \times I)$ is a collar of $(c \times \{0\}) \cup (p \times
   I) \cup (c' \times \{1\})$.

\end{prop}

\proof Consider the regions of $(I \times I) \backslash \Gamma$
abutting $I \times \{0\}$.  Since $S \cap (T \times \partial I)$
consist of curves essential in both $S$ and $(T \times \partial I)$
there must be an unlabelled such region.  We call this region $R_0$.
Similarly, there must be an unlabelled region, which we call $R_1$,
that abuts $I \times \{1\}$.  Because of the conditions on
$S \cap (T \times \partial I)$ we may assume that for $(s, t)$ in
$R_0$ or $R_1$, $S_s \cap T_t$ consists of curves essential in both
$S_s$ and $T_t$.  We are interested in monotone paths from $R_0$ to
$R_1$.  

Our argument breaks into two cases: Either there is a path in $I
\times I$ from $R_0$ to $R_1$ that avoids labelled regions or there is
no such path.

\begin{figure}[ht!]
\centerline{\epsfxsize=1.5in \epsfbox{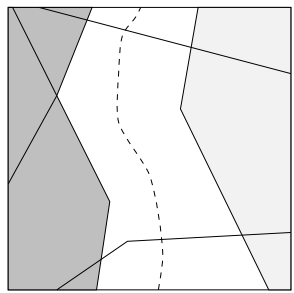}}
\caption{A monotone path through unlabelled regions}
\label{gmb.eps}
\end{figure}

{\bf Case 1}\qua There is a path beginning in $R_0$ and ending in
$R_1$ that traverses only unlabelled regions and edges between such
regions.  See Figure \ref{gmb.eps}.

It follows from Remark \ref{rem:mon} that there is a monotone path
$\alpha$ beginning in $R_0$ and ending in $R_1$ that traverses only
unlabelled regions and edges between such regions.

Consider the effect of traversing an edge $e$ from one unlabelled
region of $(I \times I) \backslash \Gamma$ to another.  Since
compression bodies do not contain essential surfaces, $S_s \cap T_t$
must contain curves essential in $S_s$.  For $(s, t)$ in an unlabelled
region, such curves must also be essential in $T_t$.  We may thus
refer to curves inessential in both $S_s$ and $T_t$ simply as
inessential and to curves essential in both $S_s$ and $T_t$ simply as
essential.

As $\alpha$ crosses $e$, the components of intersection of $S_s \cap
T_t$ change in the same way that level curves change, as we rise from
being below to being above a maximum, minimum or saddle point.  If the
point of tangency corresponding to $e$ corresponds to a maximum or
minimum, then an inessential curve appears or disappears.  If the
point of tangency corresponds to a saddle singularity, then either one
curve is banded to itself or two curves are banded together.  Note
that $S_s \cap T_t$ is separating in $T_t$, hence the essential curves
come in pairs parallel in $T_t$.  Thus if a curve is banded to itself,
then any resulting essential curves are parallel, in $T_t$, to other
essential curves.  If two essential curves are banded to each other,
then an inessential curve results, but essential curves must remain.
In both cases, the slope, in $T_t$, of the essential curves is
unaffected, as $t$ increases.

Since $\alpha$ is monotone, we obtain a function $f\co  I \rightarrow I$
by requiring $f(t)$ to be the value $s$ such that $(s, t)$ is in the
image of $\alpha$.  We may now isotope $S \cap (T \times I)$ by
isotoping $S \cap (T \times \{t\})$ to $S_{f(t)} \cap T_t$ for each $t
\in I$.  Thus $\alpha$ describes an isotopy of $S$.  We will assume in
the following that this isotopy has been performed.  The preceding
paragraph tells us how $S \cap (T \times \{t\})$ changes as $t$
increases.  We must reconstruct $S \cap (T \times I)$ from these level
sets.

As $\alpha$ traverses an unlabelled region, the curves of intersection
$S \cap T_t$ sweep out annuli.  As $\alpha$ crosses an edge
corresponding to a tangency corresponding to a maximum or minimum an
inessential curve appears or disappears.  If such a curve appears,
this indicates the appearance of a disk.  If such a curve disappears,
this indicates that an inessential curve is capped off by a disk.
Note that since there are no disks for $t =0$ or for $t = 1$, these
disks merely cap off inessential curves.  Converserly, inessential
curves are always capped off by such disks.

Suppose that $\alpha$ crosses an edge corresponding to a tangency
corresponding to a saddle point at $(s_0, t_0)$. If an
inessential curve is banded to itself, then there are two
possibilities:

\noindent
{\rm(1)}\qua  Two inessential curves result.  These curves are nested in $T_{t_0
+ \epsilon}$.  However, one of these inessential curves lies in the
subdisk of $S$ bounded by the other.  Hence the
appearance of the new inessential curve does not indicate the
appearance of another disk.  The appearance of the new inessential curve
merely indicates a saddle singularity, with respect to projection
onto $t$, in the imbedding of a single disk.  Thus we may ignore this
phenomenon.

\noindent
{\rm(2)}\qua  Two essential curves result.  Since the inessential curve bounds a
disk in $S$, the components of $S \cap T_t$
affected piece together to form an annulus that is parallel into
$T_{t_0 + \epsilon}$.

If an essential curve is banded to itself, then either one essential
curve or one inessential and one essential curve result.  But the
first possibility can't occur here, because here $S \cap T_t$ is
separating for all $t$.  Thus one inessential and one essential curve
result.  The inessential curve bounds a disk in $S$.  Hence the
appearance of the new inessential curve merely indicates a saddle
singularity, with respect to projection onto $t$, in the imbedding of
an essential annulus.  Thus we may again ignore this phenomenon.

Finally, suppose that a pair of essential curves disappears as two
essential curves are banded to each other to produce an inessential
curve. The inessential curve bounds a disk in $S$.  Thus the
components of $S \cap T_t$ affected piece together to form an annulus
that is parallel into $T_{t_0 - \epsilon}$.

Since there are no inessential curves for $t = 0$ or $t = 1$, we see
that, after further isotopies, $S \cap (T \times I)$ consists of
spanning annuli and annuli parallel into $T \times \partial I$.

\begin{figure}[ht!]
\centerline{\epsfxsize=1.5in \epsfbox{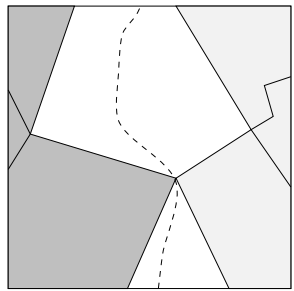}}
\caption{A monotone path through unlabelled regions}
\label{gmc.eps}
\end{figure}

{\bf Case 2}\qua There is no path beginning in $R_0$ and ending in
$R_1$ that traverses only unlabelled regions and edges between such
regions.  See Figure \ref{gmc.eps}.

In this case the regions of $(I \times I) \backslash \Gamma$ labelled
$V$ extending from the left edge of $I \times I$ must meet the regions
of $(I \times I) \backslash \Gamma$ labelled $W$ extending from the
right edge of $I \times I$.  By Remark \ref{rem:adj} this can only
happen at a vertex.  Thus Lemma \ref{lem:vertex} applies.

We consider $\tilde G$ as in C) of Lemma \ref{lem:vertex}.  Denote the
vertices by $v_1, v_2$ and the edges by $e_1, e_2, e_3, e_4$.  Suppose
first that one edge, say $e_1$ has both ends on $v_1$.  Then since
$\tilde G$ is connected, and since $v_1$ and $v_2$ have valence four,
there must be exactly two edges $e_2, e_3$, each with one end on $v_1$
and one end on $v_2$.  It follows that $e_4$ must have both ends on
$v_2$.  By Lemma \ref{lem:vertex} both $v_1 \cup e_1$ and $v_2 \cup
e_4$ are essential.  Since these curves lie on a torus, they are
parallel.  The other possibility is that all edges have one end on
$v_1$ and one end on $v_2$.

{\bf Subcase 2.1}\qua $e_1$ has both ends on $v_1$, $e_4$ has both ends on
$v_2$ and $e_2$ and $e_3$ are parallel.

\begin{figure}[ht!]
\centerline{\epsfxsize=1.5in \epsfbox{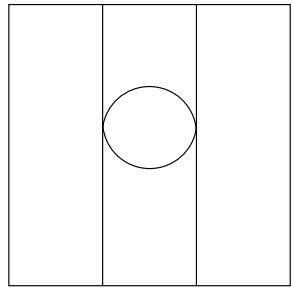}}
\caption{$\tilde G$ contains a bigon}
\label{gm2.ps}
\end{figure}

Here $\tilde G$ is as in Figure \ref{gm2.ps}.

\begin{figure}[ht!]
\centerline{\epsfxsize=\hsize \epsfbox{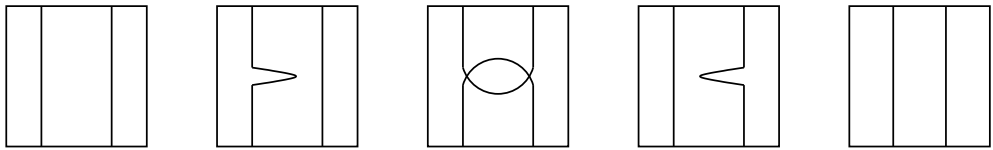}}
\caption{Adding a 1--handle}
\label{gm3.ps}
\end{figure}

This implies that as $v_0$ is crossed going from one unlabelled region
to another, $S \cap T_t$ changes as in Figure
\ref{gm3.ps}.

But this corresponds to the addition of a $1$--handle.  The core of the
$1$--handle can be isotoped into $T_{t_0}$.  Since there are essential
disks for both $V$ and $W$ cut out by $\tilde G$, $S \backslash
(T \times [t_0 - \epsilon, t_0 + \epsilon])$ is incompressible in $M
\backslash (T \times [t_0 - \epsilon, t_0 + \epsilon])$.  In
particular, we may assume that $S \cap (T \times \{t_0 - \epsilon, t_0
+ \epsilon\})$ consists of curves essential in $S$ and in
$T_{t_0 \pm \epsilon}$.  This implies that $S \cap (T \times
[t_0 - \epsilon, t_0 + \epsilon])$ is as described in option (1).

To see that option (1) holds for all of $S \cap (T \times I)$,
note that $S \cap ((T \times I) \backslash (T \times [t_0 -
\epsilon, t_0 + \epsilon]))$ consists of incompressible annuli.
Spanning annuli merely extend components of $S \cap (T \times
[t_0 - \epsilon, t_0 + \epsilon])$.  Annuli parallel into $T_{t_0 \pm
\epsilon}$ create annuli parallel into $T \times \partial I$ unless
they meet the compressible component of $S \cap (T \times [t_0 -
\epsilon, t_0 + \epsilon])$. 

\begin{figure}[ht!]
\centerline{\epsfxsize=5.0in \epsfbox{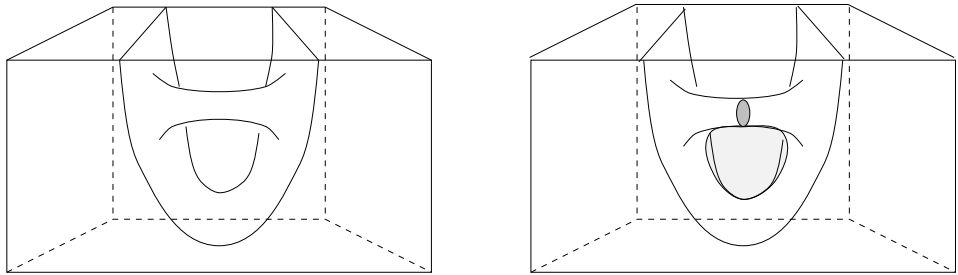}}
\caption{Compressible component meets annulus parallel into
$T_{t_0 \pm \epsilon}$}
\label{gmx.eps}
\end{figure}

\begin{figure}[ht!]
\centerline{\epsfxsize=5.0in \epsfbox{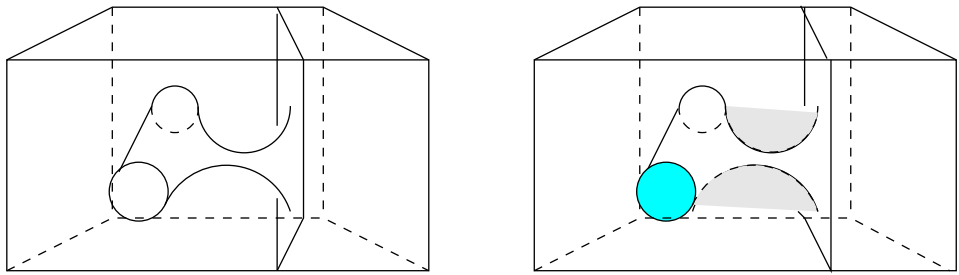}}
\caption{Another stabilized possibility}
\label{gmx1.eps}
\end{figure}

\begin{figure}[ht!]
\centerline{\epsfxsize=5.0in \epsfbox{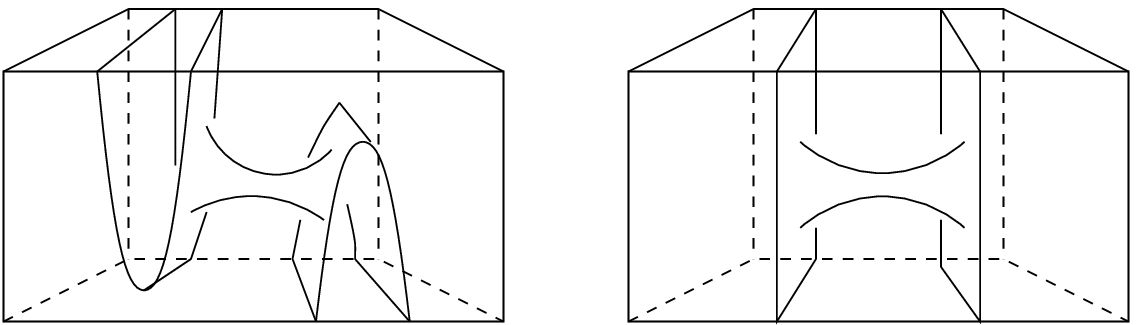}}
\caption{Two isotopic possibilities}
\label{gmx2.eps}
\end{figure}

The question is thus merely whether the ends of the compressible
component of $S \cap (T \times [t_0 - \epsilon, t_0 +
\epsilon])$ meet spanning annuli or annuli parallel into $T_{t_0 \pm
\epsilon}$.  There are a number of possibilities.  But all
possibilities are either as described in option (1) or lead to a
contradiction, by implying that $S$ is stabilized.  Some options are
pictured in Figure \ref{gmx.eps}, Figure \ref{gmx1.eps}, Figure
\ref{gmx2.eps}.  

{\bf Subcase 2.2}\qua $e_1$ has both ends on $v_1$, $e_4$ has both ends on
$v_2$ and $e_2$ and $e_3$ are not parallel.

\begin{figure}[ht!]
\centerline{\epsfxsize=1.4in \epsfbox{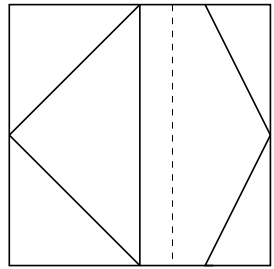}}
\caption{$\tilde G$ contains a triangle (dashed arc included along
with $\tilde G$ to indicate that a bicoloring, as required, is
possible)}
\label{gm4.ps}
\end{figure}

Here $\tilde G$ is as in Figure \ref{gm4.ps}.

\begin{figure}[ht!]
\centerline{\epsfxsize=\hsize \epsfbox{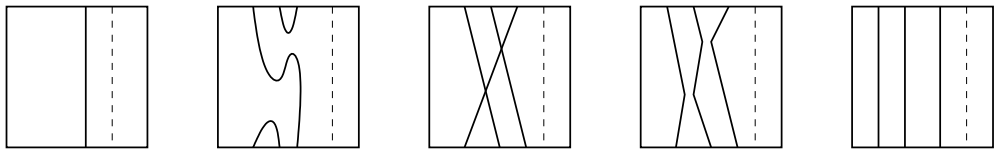}}
\caption{Adding an ``essential'' 1--handle}
\label{gm5.ps}
\end{figure}

This implies that as $v_0$ is crossed going from one unlabelled region
to another, $S \cap T_t$ changes as in Figure
\ref{gm5.ps}.

Thus $S \cap T \times [t_0 - \epsilon, t_0 + \epsilon]$ consists of
spanning annuli together with one component that is obtained as
follows: an annulus that is parallel into $T_{t_0 \pm \epsilon}$ is
tubed to a spanning annulus.  This is precisely the type of
compressible component that arises in Case 2.1 if an annulus parallel
into $T_{t_0-\epsilon}$ is attached to one end of the compressible
component and one end of a spanning annulus.  Thus as in Subcase 2.1, all
unstabilized configurations involving $S \cap T \times [t_0 -
\epsilon, t_0 + \epsilon]$ are as described in option (1).

{\bf Subcase 2.3}\qua All edges have one end on $v_1$ and the other on $v_2$.

\begin{figure}[ht!]
\centerline{\epsfxsize=1.5in \epsfbox{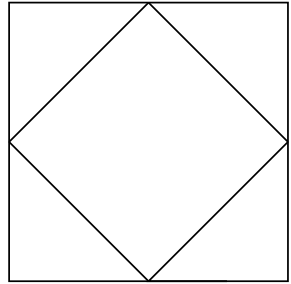}}
\caption{$\tilde G$ contains a square}
\label{gm6.ps}
\end{figure}

Note that $S$ is separating.  This induces a bicoloring of
$T_{t_0} \backslash S$.  Such a bicoloring is not possible if
two edges are parallel.  This forces $\tilde G$ to be as in Figure
\ref{gm6.ps}.

\begin{figure}[ht!]
\centerline{\epsfxsize=\hsize \epsfbox{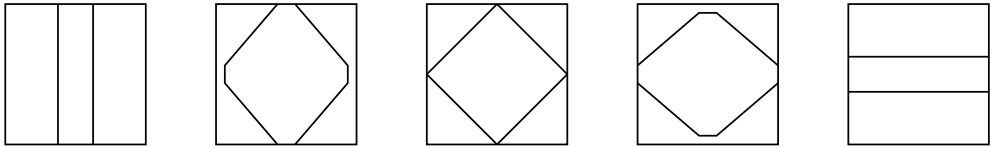}}
\caption{Another way of adding a 1--handle}
\label{gm7.ps}
\end{figure}

This implies that as $v_0$ is crossed going from one unlabelled region
to another, $S \cap T_t$ changes as in Figure
\ref{gm7.ps}.

In particular, $S \cap T \times [t_0 - \epsilon, t_0 + \epsilon]$ is
as in option (2).  Again, $S \cap ((T \times I) \backslash (T \times
[t_0 - \epsilon, t_0 + \epsilon]))$ consists of incompressible annuli.
If an end of $S \cap T \times [t_0 - \epsilon, t_0 + \epsilon]$ meets
an annulus that is not isotopic to a spanning annulus, then both ends
meet this annulus.  But this would imply that $S \cap (T \times I)$ is
disconnected.  Hence the ends of $S \cap T \times [t_0 - \epsilon, t_0
+ \epsilon]$ meet spanning annuli and $S \cap (T \times I)$ is as in
option (2).  See Figure \ref{gm7.ps}.  \endproof

\section{Fitting the pieces together}

We here prove the main theorem.  As it turns out, the hard work is
already done.  It remains to fit the results together.  

Recall that graph manifolds with empty characteristic submanifolds are
Seifert fibered spaces.  The corresponding result for Seifert fibered
spaces was obtained by Y. Moriah and the author in \cite{MS} together
with \cite{Sch1} and \cite{Sch3}.  Generalized graph manifolds with
empty characteristic submanifolds are either Seifert fibered spaces or
product manifolds of the form $(\surface) \times I$.  The corresponding
result for such manifolds was obtained by M. Scharlemann and
A. Thomspon in \cite{ST1}.

\proof[Proof of Theorem \ref{thm:premain}]
If ${\cal E} = \emptyset$, then $S$ is either pseudohorizontal
or pseudovertical, by the main theorems of \cite{ST1}, \cite{MS},
\cite{Sch1} and \cite{Sch3}.  Thus we may assume in what follows that
${\cal E} \neq \emptyset$.

 Let $\cal T$ be the collection of decomposing annuli
and tori.  Let $N$ be a component of $M \backslash \cal T$.  By Lemma
\ref{lem:C}, $S \cap N$ is incompressible unless $N$ is the active
component.  Thus if $N$ is not the active component, then $S \cap N$
is as required.  Specifically, if $N$ is a vertex manifold, $S \cap N$
is either horizontal or vertical.  And if $N$ is an edge manifold,
then $S \cap N$ consists of incompressible annuli.

If $N$ is the active component, we must consider the possibilities: If
$N$ is a vertex manifold that is Seifert fibered and has non empty
external boundary then Proposition \ref{prop:eb} applies.  If $N$ is a
vertex manifold that is a product, then Lemma \ref{lem:prod} applies.
Thus in these two cases, we may rechoose the decomposing annuli or
tori so that an edge manifold becomes the active component.
Similarly, if $N$ is vertex manifold that is a Seifert fibered space
with no external boundary, then either $S \cap N$ is pseudohorizontal,
or we may rechoose the decomposing tori so that an edge manifold
becomes the active component.

Finally, if $N$ is an edge manifold, then $S \cap N$ is as required by
Proposition \ref{prop:acs} and Proposition \ref{prop:tcs}.  
\endproof

Now Theorem \ref{thm:main} and Theorem \ref{thm:postmain} follow from
Theorem \ref{thm:premain} along with the constructions of
untelescoping and amalgamation.  Recall Theorem \ref{thm:ind} which
states that, given an irreducible Heegaard splitting, the amalgamation
of the weak reduction of this Heegaard splitting yields the original
Heegaard splitting.

Note that the amalgamation of strongly irreducible Heegaard splittings
of generalized graph manifolds interferes with the structure of the
Heegaard splittings on the generalized graph manifolds that are being
amalgamated.  For this reason, we can't prove \ref{thm:premain}
without the hypothesis of strong irreducibility.  Recall Example
\ref{ex:nostruc}, which illustrates the obstruction.  

\proof[Proof of Theorem \ref{thm:main}] Consider the irreducible Heegaard
splitting $M = V \cup_S W$.  Let $M = (V_1 \cup_{S_1} W_1) \cup_{F_1}
\dots \cup_{F_{n-1}} (V_n \cup_{S_n} W_n)$ be a weak reduction of $M =
V \cup_S W$.  Denote the decomposing tori by ${\cal T}$.  By Lemma
\ref{lem:inc}, $\cup_i F_i$ can be isotoped so that for each vertex
manifold $M_v$ of $M$, $F_i \cap M_v$ is either horizontal or vertical
and so that for each edge manifold $M_e$, $F_i \cap M_e$ consists of
incompressible tori and essential annuli.  Hence, cutting $M$ along
$\cup_i F_i$ consists of generalized totally orientable graph
manifolds.

Since $M = (V_1 \cup_{S_1} W_1) \cup_{F_1} \dots \cup_{F_{n-1}} (V_n
\cup_{S_n} W_n)$ is a weak reduction of $M = V \cup_S W$, each $M_i =
V_i \cup_{S_i} W_i$ is strongly irreducible.  Thus $M_i = V_i
\cup_{S_i} W_i$ is a strongly irreducible Heegaard splitting of a
generalized graph manifold and Theorem \ref{thm:premain} applies.  The
result now follows by juxtaposing the strongly irreducible Heegaard
splittings of the generalized graph manifolds.

The Euler characteristic calculation follows directly from the
definition of a weak reduction of a Heegaard splitting along with the
fact that here $\chi(\partial M) = 0$.  \endproof

\proof[Proof of Theorem \ref{thm:postmain}] This follows immediately from
Theorems \ref{thm:ind} and \ref{thm:main}. \endproof

Note that a horizontal (or a pseudohorizontal, respectively) surface
in a vertex manifold corresponds to a foliation of that manifold as a
surface bundle over the circle (or to a foliation of that manifold
minus a vertical solid torus as a surface bundle over the circle).  If
the splitting surface of a Heegaard splitting can be isotoped so as to
be horizontal or pseudohorizontal in two adjacent vertex manifolds,
then the glueing data must satisfy certain restrictions.  It is thus
often the case that the splitting surface of a strongly irreducible
Heegaard splitting of a graph manifold can't be isotoped to be
horizontal or pseudohorizontal in adjacent vertex manifolds.
Consequently, the canonical example of a Heegaard splitting of a graph
manifold would seem to be a Heegaard splitting obtained by taking
Heegaard splittings of the vertex manifolds with pseudovertical
splitting surfaces and amalgamating these to obtain a Heegaard
splitting of the graph manifold.

\end{document}